\documentclass[fleqn]{article}
\usepackage{amssymb,latexsym,theorem}
\newtheorem{lemma}{Lemma}[section]
\newtheorem{prop}[lemma]{Proposition}
\newtheorem{theo}[lemma]{Theorem}
\newtheorem{co}[lemma]{Corollary}

\def\pr{\noindent{\bf Proof. }}
\def\eop{\hspace*{\fill}$\Box$}
\title{Embedded polar spaces revisited}
\author{Antonio Pasini}
\date{}

\begin{document}
\maketitle

\begin{abstract}
In this paper we introduce generalized pseudo-quadratic forms and develope some theory for them. Recall that the codomain of a $(\sigma,\varepsilon)$-quadratic form is the group $\overline{K} := K/K_{\sigma,\varepsilon}$, where $K$ is the underlying division ring of the vector space on which the form is defined and $K_{\sigma,\varepsilon} := \{t-t^\sigma\varepsilon\}_{t\in K}$. Generalized pseudo-quadratic forms are defined in the same way as $(\sigma,\varepsilon)$-quadratic forms but for replacing $\overline{K}$ with a quotient $\overline{K}/\overline{R}$ for a subgroup $\overline{R}$ of $\overline{K}$ such that $\lambda^\sigma\overline{R}\lambda = \overline{R}$ for any $\lambda\in K$. In particular, every non-trivial generalized pseudo-quadratic form admits a unique sesquilinearization, characterized by the same property as the sesquilinearization of a pseudo-quadratic form. Moreover, if $q:V\rightarrow \overline{K}/\overline{R}$ is a non-trivial generalized pseudo-quadratic form and $f:V\times V\rightarrow K$ is its sesquilinarization, the points and the lines of $\mathrm{PG}(V)$ where $q$ vanishes form a subspace $S_q$ of the polar space $S_f$ associated to $f$. After a discussion of quotients and covers of generalized pseudo-quadratic forms we prove the following: let $e:S\rightarrow \mathrm{PG}(V)$ be a projective embedding of a non-degenerate polar space $S$ of rank at least $2$; then $e(S)$ is either the polar space $S_q$ associated to a generalized pseudo-quadratic form $q$ or the polar space $S_f$ associated to an alternating form $f$. By exploiting this theorem we also obtain an elementary proof of the following well known fact: an embedding $e$ as above is dominant if and only if either $e(S) = S_q$ for a pseudo-quadratic form $q$ or $\mathrm{char}(K)\neq 2$ and $e(S) = S_f$ for an alternating form $f$. 
\end{abstract}

\section{Introduction}

\subsection{Polar spaces and their embeddings}\label{Polar-Intro}

We presume that the reader is familiar with the theory of polar spaces and their projective embeddings. We refer to Tits \cite[Chapters 7 and 8]{Tits} and Buekenhout and Cohen \cite[Chapters 7-10]{BC} for this topic, but we warn the reader that there are some differences between the setting and the `philosophy' chosen by Tits \cite{Tits} and the approach of Buekenhout and Cohen \cite{BC}. To begin with, the definition of polar space adopted in \cite{BC} (which is the same as
in Buekenhout and Shult \cite{BS}) is more general than that of Tits \cite{Tits}: a polar space as defined by Tits \cite[Chapter 7]{Tits} is a non-degenerate polar space of finite rank in the sense of \cite{BC}. In this paper we shall stick to the definition of \cite{BC}, according to which a polar space is a point-line geometry $S = (P, L)$ such that for every point $p\in P$ and every line $l\in L$, the point $p$ is collinear with either all or just one of the points of $l$. The notion of projective embedding used in \cite[Chapter 8]{Tits} also looks more restrictive than that of \cite{BC}, although those two notions are in fact equivalent, as we will see in a few lines. According to \cite{BC}, an embedding of a polar space $S = (P, L)$ is an injective mapping $e$ from the point-set $P$ of $S$ to the set of points of the projective geometry $\mathrm{PG}(V)$ of a vector space $V$, such that $e$ maps every line of $S$ surjectively onto a line of $\mathrm{PG}(V)$ and $e(P)$ spans $\mathrm{PG}(V)$ (compare our definition of embeddings in Subsection \ref{Emb}), while Tits \cite{Tits} also assumes the following:

\begin{itemize}
\item[$(*)$] The image $e(S) = (e(P),e(L))$ of $S$ by $e$ is a subspace of the polar space $S_f$ associated to a reflexive sesquilinear form $f:V\times V\rightarrow K$.
\end{itemize}
Needless to say, $K$ is the underlying division ring of $V$. We warn that in $(*)$ the form $f$ is allowed to be degenerate. As for the definition of subspaces, we refer the reader to Subsection \ref{Subspaces} of this paper. 

However, as we said above, these two definitions of embedding are practically the same. Indeed:    

\begin{theo}\label{Theo 1}{\rm [Buekenhout and Cohen \cite[Chapter 9]{BC}]}
Let $e$ be a projective embedding of a polar space $S$, in the sense of {\rm \cite{BC}} (and of this paper).    
Suppose that $S$ is non-degenerate of rank at least $2$. Then $(*)$ holds for $e$. 
\end{theo}

To my knowledge, the earliest version of Theorem \ref{Theo 1} that has appeared in the literature is due to Bueknehout and Lef\'{e}vre \cite{BL}. Only polar spaces of rank 2 are considered by Buekenhout and Lef\'{e}vre \cite{BL}, but their proof also holds for higher rank polar spaces, modulo a few obvious adjustments. 

In view of the next theorem, we need a definition. Referring to Subsection \ref{Emb} for the definition of quotients and covers of embeddings, we say that a projective embedding of a polar space $S$ is {\em dominant} if it is not a proper quotient of any other embedding of $S$. In other words, it is not properly covered by any other embedding. An embedding $e$ of $S$ is {\em absolutely initial} if all projective embeddings of $S$ are quotients of $e$. (Both these definitions will be stated again in Subsection \ref{Emb}, in a more general context.) Clearly, initial embeddings are also dominant.  

\begin{theo}\label{Theo 3}{\rm [Tits \cite[8.6]{Tits}]}
Let $S$ be a non-degenerate polar space of rank at least $2$ and let $e:S\rightarrow \mathrm{PG}(V)$ be a projective embedding of $S$. Let $f$ be as in $(*)$. Then $e$ is dominant if and only if one of the following holds:
\begin{itemize}
\item[$(1)$] The form $f$ is alternating, the underlying field of $V$ has characteristic other than $2$ and $e(S) = S_f$.
\item[$(2)$] The image $e(S)$ of $S$ is the polar space $S_q$ associated to a non-singular pseudo-quadratic form $q$ such that $f$ is the 
sesquilinearization of $q$.
\end{itemize}
Moreover, if $e$ is dominant then it is also absolutely initial, except for two exceptional cases where $S$ has rank $2$. 
\end{theo}

The two exceptional cases mentioned above will be described later in this paper (Section 6, Theorem \ref{main5}). 
   
We now turn to the most important theorem of the theory of polar spaces. 

\begin{theo}\label{Theo 2}{\rm [Tits \cite[Chapter 8]{Tits}, Buekenhout and Cohen \cite[Chapter 8]{BC}]}
Let $S$ be a non-degenerate polar space of rank at least $3$. Suppose that the planes of $S$ are desarguesian. When $S$ has rank $3$ and every line of $S$ belongs to exactly two planes, suppose moreover that the planes of $S$ are Pappian. Then $S$ admits a projective embedding. 
\end{theo}

The way to prove Theorem \ref{Theo 2} is the main difference between \cite{Tits} and \cite{BC}. Tits \cite{Tits} constructs an embedding of $S$ by a free construction where vectors spaces associated to the singular subspaces of $S$ containing a given point of $S$ are amalgamated so that to obtain a vector space $\overline{V}$ which, extended by adding two copies of the underlying division ring $K$ of $S$, yields a vector space $\widetilde{V} = \overline{V}\oplus V(2,K)$ which hosts an embedding $\tilde{e}$ of $S$. The embedding $\tilde{e}$ constructed in that way is absolutely initial. Explicitly, let $\tilde{f}$ be the reflexive sequilinear form on $\widetilde{V}$ such that $\tilde{e}$ is a subspace of $S_{\tilde{f}}$ (see $(*)$). If $\tilde{e}(S) = S_{\tilde{f}}$ then $\tilde{f}$ is non-degenerate and $\tilde{e}$ is the unique projective embedding of $S$. Otherwise, $\tilde{f}$ is the sesquilinearization of a non-singular pseudo-quadratic form $\tilde{q}$, we have $\tilde{e}(S) = S_{\tilde{q}}$ and all projective embeddings of $S$ arise as quotients of $\tilde{e}$ over a subspace of the radical $\mathrm{Rad}(\tilde{f})$ of $\tilde{f}$. Thus we also have a complete classification of projective embeddings of non-degenerate polar spaces of rank at least 3.    

The proof chosen by Bukenehout and Cohen \cite{BC} is different. While Tits's proof is rather algebraic in flavour, the proof by Buekenhout and Cohen is completely geometric. Following the original approach by Veldkamp \cite{Veldkamp}, they prove that the family of hyperplanes of $S = (P, L)$ (see Subsection \ref{Subspaces} for the definition of hyperplanes) forms a projective space, say it ${\cal V}(S)$, called the {\em Veldkamp space} of $S$. The hyperplanes of $S$ are the points of ${\cal V}(S)$ while the lines of ${\cal V}(S)$ are families of hyperplanes consisting of all hyperplanes of $S$ containing the intersection of two given hyperplanes. As $S$ is non-degenerate by assumption, for every point $p\in P$ the set of points of $S$ collinear with $p$ is a hyperplane of $S$, hence a point of ${\cal V}(S)$, usually denoted by the symbol $p^\perp$. Let $\hat{e}$ be the mapping from the point-set of $S$ to the set of points of ${\cal V}(S)$ defined by setting $\hat{e}(p) = p^\perp$ for every $p\in P$. Then $\hat{e}$ is an embedding of $S$ in the subspace $\widehat{V}$ of ${\cal V}(S)$ spanned by $\hat{e}(P)$. We call $\hat{e}$ the {\em Veldkamp embedding} of $S$. 

In a sense, the Veldkamp embedding $\hat{e}$ is the counterpart of the initial embedding $\tilde{e}$ constructed by Tits. Indeed, while $\tilde{e}$ covers all embeddings of $S$, the Veldkamp embedding is covered by all of them. In short, $\hat{e}$ is terminal. 

Starting with $\hat{e}$ instead of $\tilde{e}$, in order to classify projective embeddings of polar spaces, we should describe all covers of $\hat{e}$, or at least the dominant ones, the remaining ones being obtainable as quotients of the latters. However, if we forbid ourselves to exploit the `only if' part of Tits's Theorem \ref{Theo 3} (since using that part of that theorem would imply to switching from \cite{BC} to \cite{Tits}), all we can say in general on $\hat{e}$ and its covers is what Theorem \ref{Theo 1} tells us. According to that theorem, if $e$ is an embedding of $S$ then $e(S)$ is a subspace of $S_f$ for a suitable reflexive sesquilinear form $f$, but it can happen that $e(S)$ is a proper subspace of $S_f$ as well as a proper overspace of $S_q$ for every pseudo-quadratic form $q$ admitting $f$ as the sequilinearization. As a consequence, if $\hat{f}$ is the $(\sigma,\varepsilon)$-sesquilinear form on $\widehat{V}$ such that $\hat{e}(S)$ is a subspace of $S_{\hat{f}}$ (see $(*)$) and $\tilde{e}$ is the initial embedding of $S$, when $\hat{e}(S)\subset S_{\hat{f}}$ we can only say that $\tilde{e}(S) = S_{\tilde{q}}$ for a suitable $(\sigma,\varepsilon)$-quadratic form $\tilde{q}$ defined on a suitable subspace $\widetilde{V}$ of $\overline{V}\oplus \overline{K}^{\sigma,\varepsilon}/K_{\sigma,\varepsilon}$ (compare Buekenhout and Cohen \cite[Theorem 10.12.5]{BC}), but we would get in troubles if asked for a more precise description of $\widetilde{V}$ valid in general, although we can give such a description in many particular cases.         

\subsection{Purposes and main result of this paper} 

The purpose of this paper is to overcame the difficulties discussed at the end of the previous subsection. We will succeed by introducing generalized pseudo-quadratic forms. 

 We recall that the codomain of a $(\sigma,\varepsilon)$-quadratic form is the group $\overline{K} := K/K_{\sigma,\varepsilon}$, where $K$ is the underlying division ring of the vector space $V$ on which the form is defined and $K_{\sigma,\varepsilon} := \{t-t^\sigma\varepsilon\}_{t\in K}$. {\em Generalized pseudo-quadratic forms}, to be introduced and discussed in Section 3, are defined in the same way as $(\sigma,\varepsilon)$-quadratic forms but for replacing $\overline{K}$ with a quotient $\overline{K}/\overline{R}$ for a subgroup $\overline{R}$ of $\overline{K}$ such that $\lambda^\sigma\overline{R}\lambda = \overline{R}$ for any $\lambda\in K$. In particular, every generalized pseudo-quadratic form admits a sesquilinearization, characterized by the same property as the sesquilinearization of a pseudo-quadratic form. As we shall prove in Section 3, the sesquilinearization of a non-trivial generalized pseudo-quadratic form is unique. Let $q:V\rightarrow \overline{K}/\overline{R}$ be a non-trivial generalized pseudo-quadratic form and let $f:V\times V\rightarrow K$ be its sesquilinarization. Then the points and the lines of $\mathrm{PG}(V)$ where $q$ vanishes form a subspace $S_q$ of $S_f$ (see Section 3). In Section 5 (Theorems \ref{main1} and \ref{qE4}) we shall obtain the following improvement of Theorem \ref{Theo 1}: 

\begin{theo}\label{Theo 4}
Let $e:S\rightarrow \mathrm{PG}(V)$ be a projective embedding of a non-degenerate polar space $S$ of rank at least $2$. Then $e(S)$ is either the polar space $S_q$ associated to a non-trivial generalized pseudo-quadratic form $q$ or the polar space $S_f$ associated to a non-degenerate alternating form $f$.
\end{theo}

We recall that the {\em hull} of an embedding $e$ is the unique dominant embedding that covers $e$ (see Subsection \ref{Emb}), uniqueness being understood modulo isomorphisms. With $e$ and $S$ as in Theorem \ref{Theo 4}, the hull of $e$ is the initial embedding of $S$, with the only exception of the two cases of rank $2$ mentioned in the last claim of Theorem \ref{Theo 3}. 

Let $e(S) = S_f$ for an alternating form $f$ and let $\tilde{e}$ be the hull of $e$. It is well known that in this case either $\tilde{e} = e$ (when $\mathrm{char}(K) \neq 2$) or $\mathrm{char}(K) = 2$ and $\tilde{e}(S) = S_{\tilde{q}}$ for a non-singular quadratic form $\tilde{q}:\widetilde{V}\rightarrow K$, where $\tilde{V} = V\oplus K$, the field $K$ being regarded as a vector field over itself with scalar multiplication $\circ:K\times K\rightarrow K$ defined as follows: $t\circ \lambda = t\lambda^2$ for every vector $t\in K$ and every scalar $\lambda\in K$.

On the other hand, let $e(S) = S_q$ for a generalized pseudo-quadratic form $q:V\rightarrow \overline{K}/\overline{R}$. Let $\circ:\overline{R}\times K\rightarrow K$ be defined as follows: $r\circ \lambda = \lambda^\sigma r\lambda$ for every $r\in \overline{R}$ and every scalar $\lambda\in K$. We will prove in Section 3 that the group $\overline{R}$ equipped with $\circ$ as the scalar multiplication is a $K$-vector space. (This amounts to say that $\overline{R}\subseteq K^{\sigma,\varepsilon}/K_{\sigma,\varepsilon}$.) Hence we can form a direct sum of $K$-vector spaces $\widetilde{V} = V\oplus \overline{R}$ and, if $f$ is the sesquilinearization of $q$, we can define a reflexive sesquiliner form $\tilde{f}:\widetilde{V}\times\widetilde{V}\rightarrow K$ by declaring that $\overline{R}\subseteq \mathrm{Rad}(\tilde{f})$ and $\tilde{f}$ induces $f$ on $V\times V$. As we shall prove in Section 4, a pseudo-quadratic form $\tilde{q}:\widetilde{V}\rightarrow\overline{K}$ can be defined on $\widetilde{V}$ admitting $\tilde{f}$ as it sesquilinearization and such that the projection $\pi:\widetilde{V}\rightarrow \widetilde{V}/\overline{R} = V$ induces an isomorphism $\pi_S$ from $S_{\tilde{q}}$ to $S_q$. So, the mapping $\tilde{e}:= \pi_S^{-1}\cdot e$ is a projective embedding of $S$ and $\pi$ is a morphism from $\tilde{e}$ to $\varepsilon$. Moreover, $\tilde{e}$ is dominant by Theorem \ref{Theo 3}, since $\tilde{e}(S) = S_{\tilde{q}}$ and $\tilde{q}$ is pseudo-quadratic. Therefore:

\begin{theo}\label{Theo 5}
The hull of $e$ is the embedding $\tilde{e}$ defined as above.
\end{theo}
{\bf Organization of the paper.} In the rest of Section 1 we recall some basics on subspaces and embeddings of point-line geometries. In Section 2 we give a summary of the theory of reflexive sesquilinear forms, pseudoquadratic forms and related polar spaces. In Section 3 we introduce generalized pseudo-quadratic forms and develope some theory for them. Quotients and covers of generalized pseudo-quadratic forms are discussed in Section 4. Section 5 is devoted to the proof of Theorem \ref{Theo 4}. Finally, in Section 6 we revisit Theorem \ref{Theo 3}. 

\subsection{Subspaces and embeddings of point-line geometries}

In this subsection we fix some terminology for point-line geometries, focusing on subspaces and projective embeddings. 

Throughout this subsection $G = (P, L)$ is a point-line geometry, with $P$ and $L$ as the point-set and the line-set respectively.  We regard lines as subsets of $P$ and we assume that no two distinct lines meet in more than one point and every line has at least two points. The {\em collinearity graph} of $G$ is the graph with $P$ as the vertex-set where two points $a, b\in P$ are declared to be adjacent when they are joined by a line of $G$. The geometry $G$ is said to be {\em connected} if its collinearity graph is connected. 

Given two point-line geometries $G = (P, L)$ and $G' = (P', L')$, an {\em isomorphism} from $G$ to $G'$ is a bijective mapping $e:P\rightarrow P'$ such that $\{e(l)\}_{l\in L} = L'$, where for a line $l\in L$ we put $e(l) := \{e(p)\}_{p\in l}$. 

\subsubsection{Subgeometries and subspaces}\label{Subspaces}

A point-line geometry $G' = (P', L')$ is a {\em subgeometry} of $G = (P, L)$ if $P'\subseteq P$ and for every line $l'\in L'$ there exists a (necessarily unique) line $l\in L$ such that $l' = l\cap P'$. If every line of $G'$ is also a line of $G$ then $G'$ is called a {\em full} subgeometry of $G$. On the other hand, if $L' = \{l\cap P' ~|~ l\in L, |l\cap P'| \geq 2\}$ then $G'$ is called the subgeometry {\em induced} by $G$ on $P'$.

A subset $P'\subseteq P$ is called a {\em subspace} of $G$ if every line of $G$ either is contained in $P'$ or meets $P'$ in at most one point. We say that a geometry $G' = (P', L')$ is a {\em subspace} of $(G, L)$ if $P'$ is a subspace of $G$ in the previous sense and $G'$ is the subgeometry induced by $G$ on $P'$. Clearly, subspaces in the latter sense are full subgeometries. 

We have mentioned hyperplanes in Subsection \ref{Polar-Intro}. A {\em hyperplane} of a point-line geometry $G = (P,L)$ is a proper subspace $H\subset P$ such that every line of $G$ either meets $H$ in a single point or it is fully contained in $H$. 

\subsubsection{Notation for vector spaces and projective spaces}\label{Emb-pre}

In view of the next subsection, it is convenient to fix some notation for vector spaces and related projective spaces. Given a vector space $V$, we denote by $\mathrm{PG}(V)$ the projective space of $1$-and $2$-dimensional vector subspaces of $V$. For a vector $v\in V-\{0\}$, we denote by $[v]$ the projective point of $\mathrm{PG}(V)$ represented by $v$. If $X$ is a subspace of $V$ we put $[X] = \{[x]\}_{x\in X-\{0\}}$, namely $[X]$ is the subspace of $\mathrm{PG}(V)$ corresponding to $X$. Given a semilinear mapping $f:V\rightarrow V'$, let $\mathrm{Ker}(f) := f^{-1}(0)$ be the kernel of $f$. We denote by $\mathrm{PG}(f)$ the mapping induced by $f$ from $\mathrm{PG}(V)-[\mathrm{Ker}(f)]$ to $\mathrm{PG}(V')$. 

\subsubsection{Projective embeddings}\label{Emb} 

Let $G = (P,L)$ be a connected point-line geometry. A {\em projective embedding} of $G$ (also called just {\em embedding} for short) is an isomorphism $e$ from $G$ to a full subgeometry $e(G) = (e(P), e(L))$ of the projective space $\mathrm{PG}(V)$ of a vector space $V$, such that $e(P)$ spans $\mathrm{PG}(V)$. 

We write $e:G\rightarrow \mathrm{PG}(V)$ to mean that $e$ is a projective embedding of $G$ in $\mathrm{PG}(V)$. If $K$ is the underlying division ring of $V$ then we say that $e$ is {\em defined over} $K$, also that $e$ is a $K$-{\em embedding}, for short. If all projective embeddings of $G$ are defined over the same division ring $K$ then we say that $G$ is {\em defined over} $K$ and we call $K$ the {\em underlying division ring} of $G$. 

Given two $K$-embeddings $e:G\rightarrow \mathrm{PG}(V)$ and $e'\rightarrow\mathrm{PG}(V')$, a {\em morphism} $f:e\rightarrow e'$ is a semilinear mapping $f:V\rightarrow V'$ such that $\mathrm{PG}(f)\cdot e = e'$. As $e'(P)$ spans $\mathrm{PG}(V')$, the mapping $f$ is surjective. If $f$ is bijective then we say that $f$ is an {\em isomorphism} from $e$ to $e'$. If a morphism $f:e\rightarrow e'$ exists then we say that $e'$ is a {\em morphic image} of $e$ (also that $e$ {\em covers} $e'$) and we write $e\geq e'$. If moreover $f$ is bijective then we write $e\cong e'$ and we say that $e$ and $e'$ are {\em isomorphic}. If the morphism $f$ is not an isomorphism then we call $f$ a {\em proper morphism} and we write $e > e'$. Note that, as $G$ is connected by assumption, if $e\geq e'$ then the morphism $f:e\rightarrow e'$ is unique up to isomorphims.  

Let $U$ be a subspace of $V$ such that $e(P)\cap[U]= \emptyset$ and $l\cap[U] = \emptyset$ for any line $l$ of $\mathrm{PG}(V)$ such that $|l\cap e(P)| \geq 2$. Let $\pi_U$ the projection of $V$ onto $V/U$. Then the mapping $e_U := \mathrm{PG}(\pi_U)\circ e$ is an embedding of $G$ in $\mathrm{PG}(V/U)$ and $\pi_U$ is a morphism from $e$ to $e_U$. We say that $U$ {\em defines a quotient} of $e$ a we call $e_U$ the {\em quotient} of $e$ over $U$. 

Clearly, if $f:e\rightarrow e'$ is a morphism then $\mathrm{Ker}(f)$ defines a quotient of $e$ and we have $e' \cong e_U$. By a little abuse, we say that $e'$ is a {\em quotient} of $e$, thus taking the word `quotient' as a synonym of `morphic image'.  

Following Tits \cite[Chapter 8]{Tits} we say that a projective embedding of $G$ is {\em dominant} if it cannot be obtained as a proper quotient from any other projective embedding of $G$.  If all $K$-embeddings of $G$ are quotient of a given $K$-embedding $e$ then we say that $e$ is $K$-{\em initial}. If moreover all embeddings of $G$ are quotients of $e$ then $e$ is said to be {\em absolutely initial}. 

Clearly, the $K$-initial embedding, if it exists, is uniquely determined up to isomorphisms. It can be characterized as the unique dominant $K$-embedding of $G$.  It is also clear that $G$ admits the absolutely initial embedding if and only if it is defined over some division ring $K$ and admits the $K$-initial embedding. 

Finally, every embedding $e$ of $G$ admits a {\em hull} $\tilde{e}$, uniquely determined up to isomorphism by the following property: $\tilde{e}\geq e'$ for every embedding $e'$ of $G$ such that $e'\geq e$. We refer the reader to Ronan \cite{Ron} for an explicit construction of $\tilde{e}$. Clearly, the hull $\tilde{e}$ of $e$ is dominant. Up to isomorphisms, it is the unique dominant embedding in the class of the embeddings that cover $e$. So, if $G$ admits the $K$-initial embedding and $e$ is defined over $K$, then $\tilde{e}$ is also $K$-initial. 

The terminology adopted in the previous definitions is essentially the same as in Tits \cite{Tits}, but we warn the reader that different terminologies are also used in the literature. For instance, dominant embeddings and $K$-initial embeddings are often called {\em relatively universal} and {\em absolutely universal} respectively (compare Kasikova and Shult \cite{KS}).   

\section{Preliminaries}

In this section we fix some notation to be used troughout this paper and recall a few basics on sesquilinear and pseudo-quadratic forms, taken from Tits \cite[Chapter 8]{Tits} and Buekenhout and Cohen \cite[Chapters 7 and 10]{BC}. In Section 3, all properties of pseudo-quadratic forms to be recalled in the present section will be rephrased in the setting of generalized pseudo-quadradic forms.    

\subsection{Admissible pairs}\label{section 2.1}

Throughout this paper $K$ is a possibly non-commutative division ring, $\sigma$ is an anti-automorphism of $K$ and $\varepsilon\in K$ is such that $\varepsilon^\sigma\varepsilon = 1$ and $t^{\sigma^2} = \varepsilon t\varepsilon^{-1}$ for any $t\in K$. Following Buekenhout and Cohen \cite[Chapter 10]{BC} we call $(\sigma, \varepsilon)$ an {\em admissible pair} of $K$. As in Tits \cite[Chapter 8]{Tits}, we set
\[K_{\sigma,\varepsilon} := \{t - t^\sigma\varepsilon\}_{t\in K}, \hspace{10 mm} K^{\sigma,\varepsilon} = \{t\in K ~|~ t = -t^\sigma\varepsilon\}. \]
Clearly $K_{\sigma,\varepsilon}$ and $K^{\sigma,\varepsilon}$ are subgroups of the additive group of $K$. Moreover
\begin{equation}\label{sigma-epsilon1}
\lambda^\sigma K_{\sigma,\varepsilon}\lambda = K_{\sigma,\varepsilon} ~~ \mbox{and}~~ \lambda^\sigma K^{\sigma,\varepsilon}\lambda = K^{\sigma,\varepsilon} ~~\mbox{ for every}~\lambda\in K-\{0\},
\end{equation}
\begin{equation}\label{sigma-epsilon2}
K_{\sigma,\varepsilon} \subseteq K^{\sigma,\varepsilon},
\end{equation} 
\begin{equation}\label{sigma-epsilon3}
\left.\begin{array}{rcl}
K^{\sigma,\varepsilon} = K & \mbox{if and only if} & \sigma = \mathrm{id}_K ~ \mbox{and} ~ \varepsilon = -1,\\
K_{\sigma,\varepsilon} = K & \mbox{if and only if} &  \sigma = \mathrm{id}_K, ~\varepsilon = -1 ~ \mbox{and} ~ \mathrm{char}(K) \neq 2. 
\end{array}\right\}
\end{equation}
The quotient group of the additive group of $K$ over $K_{\sigma,\varepsilon}$ is denoted by $K^{(\sigma,\varepsilon)}$ in \cite{Tits}. In this paper we shall denote it by the symbol $\overline{K}$: 
\begin{equation}\label{Kse2}
\overline{K} :=  K^{(\sigma,\varepsilon)} = K/K_{\sigma,\varepsilon}.
\end{equation}
We will also adopt the following convention. Given $t\in K$ we denote by $\bar{t}$ the element of $\overline{K}$ represented by $t$:
\[\bar{t}:= t + K_{\sigma,\varepsilon}.\]
Accordingly, $\overline{t+s} = t+s + K_{\sigma,\varepsilon}$, $\overline{ts} = ts + K_{\sigma,\varepsilon}$ and $\bar{0}$ is the null element of $\overline{K}$.

\subsubsection{Pairs of trace type} 

Clearly, if $(\sigma,\varepsilon)$ is an admissible pair of a division ring $K$ then the pair $(\sigma,-\varepsilon)$ is also admissible. So, we can consider the groups $K_{\sigma,-\varepsilon} = \{t+t^\sigma\varepsilon\}_{t\in K}$ and $K^{\sigma,-\varepsilon} = \{t \in K ~|~ t = t^\sigma\varepsilon\}$. According to (\ref{sigma-epsilon2}),  $K_{\sigma,-\varepsilon}\subseteq K^{\sigma,-\varepsilon}$. Following Buekenhout and Cohen \cite{BC}, when $K_{\sigma,-\varepsilon} = K^{\sigma,-\varepsilon}$ we say that the pair $(\sigma,\varepsilon)$ is of {\em trace type}.  

The following is well known (see Tits \cite[Chapter 8]{Tits}, also Buekenhout and Cohen \cite[Chapter 10]{BC}). 

\begin{lemma}\label{trace type pair}
Assume that either $\mathrm{char}(K) \neq 2$ or $\mathrm{char}(K) = 2$ but $\sigma$ acts non-trivially on the center $Z(K)$ of $K$. Then, for every element $\varepsilon\in K$ forming an admissible pair with $\sigma$, the pair $(\sigma,\varepsilon)$ is of trace type. 
\end{lemma}

\subsubsection{A scalar multiplication in the group $\overline{K}$}\label{scalar-mult-sec}

According to (\ref{sigma-epsilon1}),  $\lambda^\sigma K_{\sigma,\varepsilon}\lambda \subseteq K_{\sigma,\varepsilon}$ for every $\lambda\in K$. Hence can define a scalar multiplication $\circ : \overline{K}\times K \rightarrow \overline{K}$ as follows: $(t+K_{\sigma,\varepsilon})\circ\lambda = \lambda^\sigma(t+K_{\sigma,\varepsilon})\lambda = \lambda^\sigma t\lambda + K_{\sigma,\varepsilon}$, namely
\begin{equation}\label{circ}
\bar{t}\circ\lambda = \overline{\lambda^\sigma t\lambda} ~~ \mbox{for any}~ \bar{t}\in \overline{K}~ \mbox{and}~ \lambda\in K.
\end{equation}
Clearly the following hold for any $\bar{t},\bar{s}\in\overline{K}$ and $\lambda, \mu\in K$:
\begin{equation}\label{property circ}
(\bar{t}\circ\lambda)\circ\mu = \bar{t}\circ(\lambda\mu) ~~ \mbox{and} ~~ (\overline{t+s})\circ\lambda = \bar{t}\circ\lambda + \bar{s}\circ\lambda.
\end{equation}
Given an element $\bar{t}\in \overline{K}$ (a subset $\overline{H}\subseteq \overline{K}$) we put $\bar{t}\circ K := \{\bar{t}\circ\lambda\}_{\lambda\in K}$ (respectively $\overline{H}\circ K := \cup_{\bar{t}\in\overline{H}}\bar{t}\circ K$). We say that an element $\bar{t}\in\overline{K}$ is a $\circ$-{\em vector} if 
\begin{equation}\label{def1}
\bar{t}\circ(\lambda+\mu) = \bar{t}\circ\lambda + \bar{t}\circ\mu ~~\mbox{for any}~ \lambda,\mu\in K.
\end{equation}
We denote by $\overline{K}^\circ$ the set of $\circ$-vectors of $\overline{K}$. It is easy to see that $\overline{K}^\circ+\overline{K}^\circ \subseteq \overline{K}^\circ$ and $\overline{K}^\circ\circ K\subseteq \overline{K}^\circ$. Moreover, $\bar{0}\in \overline{K}^\circ$ and $-\overline{K}^\circ = \overline{K}^\circ$. Thus, $\overline{K}^\circ$ can be regarded as a right $K$-vector space, with $\circ$ taken as the scalar multiplication. 

The next lemma is essentially the same as Lemma 10.2.2 of Buekenhout and Cohen \cite{BC}. We leave the proof for the reader. 

\begin{lemma}\label{lemma1}
We have $\overline{K}^\circ = K^{\sigma,\varepsilon}/K_{\sigma,\varepsilon}$.
\end{lemma}


The next corollary is well known. (See Tits \cite[Chapter 8]{Tits}, also Bukenhout and Cohen \cite[Chapter 10]{BC}, and recall that $\overline{K}^\circ = K^{\sigma,\varepsilon}/K_{\sigma,\varepsilon}$ by Lemma \ref{lemma1}.)

\begin{co}\label{cor1}
Both the following hold.

\begin{itemize}
\item[$(1)$] $\overline{K}^\circ = \{\bar{0}\}$ if and only if the pair $(\sigma,\varepsilon)$ is of trace type. 
\item[$(2)$] $\overline{K}^\circ = \overline{K}$ if and only if $K^{\sigma,\varepsilon} = K$.
\end{itemize}
\end{co}

\subsubsection{Closed subgroups of $\overline{K}$} 

We say that a subgroup $\overline{H}$ of $\overline{K}$ is {\em closed} with respect to the scalar multiplication $\circ$ (also $\circ$-{\em closed} or just {\em closed}, for short) if $\overline{H}\circ K \subseteq \overline{H}$. 

Clearly $\overline{K}$, the vector space $\overline{K}^\circ$ and all of its subspaces are closed subgroup of $\overline{K}$. We are not going to discuss properties of closed subgroups here. We only mention the following, to be exploited in Section 3.    

Let $\overline{H}$ be a closed subgroup of $\overline{K}$. The scalar multiplication $\circ$ of $\overline{K}$ naturally induces a scalar multiplication on the quotient group $\overline{K}/\overline{H}$, which we shall denote by the same symbol $\circ$ used for the scalar multiplication of $\overline{K}$. Explicitly, 
\begin{equation}\label{circ-quot}
(\bar{t}+\overline{H})\circ\lambda := \bar{t}\circ\lambda + \overline{H} ~~ \mbox{for every}~\bar{t}\in \overline{K}.
\end{equation}
It is easy to see that this definition is consistent, namely the coset $\bar{t}\circ\lambda + \overline{H}$ does not depend on the choice of the representative $\bar{t}$ of $\bar{t}+\overline{H}$. Moreover, the scalar multiplication defined on $\overline{K}/\overline{H}$ in this way satisfies identities similar to (\ref{property circ}).

\subsubsection{Proportionality of admissible pairs}\label{proportional pairs sec}

Given an admissible pair $(\sigma,\varepsilon)$ of $K$ and a nonzero scalar $\kappa\in K-\{0\}$, let $\varepsilon' := \kappa\kappa^{-\sigma}\varepsilon$ an let $\sigma'$ be the anti-automorphism of  $K$ defined as follows:
\[t^{\sigma'} := \kappa t^\sigma \kappa^{-1} ~~~ \mbox{for every} ~ t\in K.\]
All claims gathered in the next lemma are well known (see Tits \cite[Chapter 8]{Tits}):

\begin{lemma}\label{proportional pairs lemma}
The pair $(\sigma',\varepsilon')$ is admissible. Moreover:

\begin{itemize}
\item[$(1)$] We have $\kappa K_{\sigma,\varepsilon} = K_{\sigma',\varepsilon'}$ and $\kappa K^{\sigma,\varepsilon} = K^{\sigma',\varepsilon'}$.
\item[$(2)$] $\kappa \lambda^\sigma t\lambda = \lambda^{\sigma'}\kappa t\lambda$ for any $t\in K$.
\end{itemize}
\end{lemma}

By $(1)$ of Lemma \ref{proportional pairs lemma}, left multiplication by $\kappa$ induces a group isomorphism from $K/K_{\sigma,\varepsilon}$ to $K/ K_{\sigma',\varepsilon'}$ as well as from $K^{\sigma,\varepsilon}/K_{\sigma,\varepsilon}$ to $K^{\sigma',\varepsilon'}/ K_{\sigma',\varepsilon'}$.

When dealing with two pairs $(\sigma,\varepsilon)$ and $(\sigma',\varepsilon')$ as above it is convenient to keep a record of them in our notation. So we put $\overline{K}^{\sigma,\varepsilon} =K/ K_{\sigma,\varepsilon}$, $\overline{K}^{\sigma',\varepsilon'} =K/ K_{\sigma',\varepsilon'}$, $\overline{K}^{\circ, \sigma,\varepsilon} =K^{\sigma,\varepsilon}/ K_{\sigma,\varepsilon}$, $\overline{K}^{\circ,\sigma',\varepsilon'} =K^{\sigma',\varepsilon'}/ K_{\sigma',\varepsilon'}$, $\bar{t}^{\sigma,\varepsilon} = t + K_{\sigma,\varepsilon}$, $\bar{t}^{\sigma',\varepsilon'} = t + K_{\sigma',\varepsilon'}$ and we denote the scalar multiplications of $\overline{K}^{\sigma,\varepsilon}$ and $\overline{K}^{\sigma',\varepsilon'}$ by the symbols $\circ_{\sigma}$ and $\circ_{\sigma'}$ respectively. This notation is admittedly rather clumsy. We will avoid it as far as possible, but in the present context we need it.

With the above notation, claim $(2)$ of Lemma \ref{proportional pairs lemma} can be rewritten as follows:
\[\kappa(\bar{t}^{\sigma,\varepsilon}\circ_\sigma \lambda) = (\kappa(\bar{t}^{\sigma,\varepsilon}))\circ_{\sigma'}\lambda =
(\overline{(\kappa t)
}^{\sigma',\varepsilon'})\circ_{\sigma'}\lambda.\]
Thus, left multiplication by $\kappa$ is an isomorphism of $K$-vector spaces from $\overline{K}^{\circ, \sigma,\varepsilon}$ to $\overline{K}^{\circ,\sigma',\varepsilon'}$. 

With $\kappa$, $(\sigma,\varepsilon)$ and $(\sigma',\varepsilon')$ as above, we write $(\sigma',\varepsilon') = \kappa \cdot (\sigma,\varepsilon)$ and we say that the pairs $(\sigma,\varepsilon)$ and $(\sigma',\varepsilon')$ are {\em proportional}. 

Clearly, if $(\sigma',\varepsilon') = \kappa\cdot(\sigma,\varepsilon)$ then $(\sigma,\varepsilon) = \kappa^{-1}\cdot (\sigma',\varepsilon')$. If moreover $(\sigma'',\varepsilon'') = \kappa'\cdot(\sigma,\varepsilon)$ then $(\sigma'',\varepsilon'') = (\kappa'\kappa)\cdot(\sigma,\varepsilon)$. It is also clear that $\kappa\cdot(\sigma,\varepsilon) = (\sigma,\varepsilon)$ if and only if $\kappa\in Z(K)$ and $\kappa^\sigma = \kappa$. 

\subsection{Reflexive sesquilinear forms}

Given a division ring $K$, a left $K$-vector space $V$ and an antiautomorphism $\sigma$ of $K$, a $\sigma$-{\em sesquilinear form} is a mapping $f:V\times V\rightarrow K$ such that
\begin{equation}\label{sesqui1}
\begin{array}{l}
f(x_1\lambda_1 + x_2\lambda_2, y_1\mu_1+ y_2\mu_2) = \\
= \lambda_1^\sigma f(x_,y_1)\mu_1 + \lambda_1^\sigma f(x_1, y_2)\mu_2 + \lambda_2^\sigma  f(x_2,y_1)\mu_1 + \lambda_2^\sigma f(x_2,y_2)\mu_2
\end{array}
\end{equation}
for any $x_1, x_2, y_1, y_2\in V$ and $\lambda_1, \lambda_2, \mu_1, \mu_2\in K$. We say that $f$ is {\em trivial} when $f(x,y) = 0$ for any choice of $x, y \in V$. Obviously, if $f$ is non-trivial then $\sigma$ is uniquely determined by (\ref{sesqui1}).   

A sesquilinear form $f$ is said to be {\em reflexive} if, for any choice of $x, y \in V$, we have $f(x,y) = 0$ if and only if $f(y,x) = 0$. It is well known (Tits \cite[Chapter 8]{Tits}) that a non-trivial $\sigma$-sesquilinear form is reflexive if and only if there exists an element $\varepsilon\in K$ such that 
\begin{equation}\label{sesqui2}
f(y,x) = f(x,y)^\sigma\varepsilon ~~ \mbox{for any choice of} ~ x, y \in V.
\end{equation}
If this is the case then $(\sigma,\varepsilon)$ is an admissible pair and $f$ is called a $(\sigma,\varepsilon)$-{\em sesquilinear form}. Clearly, the element $\varepsilon$ satisfying (\ref{sesqui2}) is unique.    

A {\em bilinear form} is a $\sigma$-sesquilinear form with $\sigma = \mathrm{id}_K$ (whence $K$ is a field, namely it is commutative). A {\em symmetric bilinear} form is an $(\mathrm{id}_K,1)$-sesquilinear form. A bilinear form $f$ is said to be {\em alternating} if
\begin{equation}\label{alt}
f(x,x) = 0 ~~\mbox{for any}~ x\in V.
\end{equation}
Non-trivial alternating forms are $(\mathrm{id}_K,-1)$-sesquilinear. Conversely, if $K$ is a field of characteristic $\mathrm{char}(K) \neq 2$ then all $(\mathrm{id}_K,-1)$-sesquilinear forms are alternating. On the other hand, let $\mathrm{char}(K) = 2$. Then $1 = -1$. In this case a $(\mathrm{id}_K,-1)$-sesquilinear form is just a symmetric bilinear form. Obviously, not all symmetric bilinear forms satisfy (\ref{alt}).   

\subsubsection{Orthogonality}\label{perp}

Given a $(\sigma,\varepsilon)$-sesquilinear form $f:V\times V\rightarrow K$, we say that two vectors $x, y \in K$ are {\em orthogonal} (with respect to $f$) if $f(x,y) = 0$. If $x$ and $y$ are orthogonal then we write $x\perp y$. Given a vector $x\in V$ we put $x^\perp := \{y \in V ~|~ y\perp x\}$ and, for a subset $X\subseteq V$, we set $X^\perp := \bigcap_{x\in X}x^\perp$. Clearly $x^\perp$ is either a hyperplane or the whole of $V$. Hence $X^\perp$ is a subspace of $V$, for any $X\subseteq V$. Note also that $\langle X\rangle^\perp = X^\perp$. We set 
\[\mathrm{Rad}(f) := V^\perp = \{x\in V~|~ x^\perp = V\}\]
and we call $\mathrm{Rad}(f)$ the {\em radical} of $f$. We say that $f$ is {\em degenerate} if $\mathrm{Rad}(f) \neq \{0\}$.  

A vector $x\in V$ is said to be {\em isotropic} for $f$ (also $f$-{\em isotropic}) if $f(x,x) = 0$, namely $x\in x^\perp$. A subset $X\subseteq V$ is {\em totally isotropic} for $f$ ({\em totally $f$-isotropic}) if $X\subseteq X^\perp$. 

Clearly, all vectors of $\mathrm{Rad}(f)$ are isotropic. We say that the form $f$ is {\em strictly isotropic} if it admits at least one isotropic vector $x\not\in \mathrm{Rad}(f)$. 

\subsubsection{Trace-valued forms}

Let $f:V\times V\rightarrow K$ be a $(\sigma,\varepsilon)$-sesquilinear form. By (\ref{sesqui2}), $f(x,x)\in K^{\sigma,-\varepsilon}$ for every $x\in V$. The form $f$ is said to be {\em trace-valued} if $f(x,x)\in K_{\sigma,-\varepsilon}$ for every $x\in V$. Two well known characterizations of trace-valued forms are gathered in the next proposition (see Tits \cite[Chapter 8]{Tits}, also Buekenhout and Cohen \cite[Chapter 10]{BC}). 

\begin{prop}\label{trace-prop2}
Let $f:V\times V\rightarrow K$ be a $(\sigma,\varepsilon)$-sesquilinear form. Then:

\begin{itemize}
\item[$(1)$] The form $f$ is trace-valued if and only if there exists a $\sigma$-sesquilinear form $g:V\times V\rightarrow K$ such that $f(x,y) = g(x,y) + g(y,x)^\sigma\varepsilon$ for all $x, y \in V$.
\item[$(2)$] Asume that $f$ is strictly isotropic. Then $f$ is trace-valued if and only if $V$ is spanned by the set of $f$-isotropic vectors.
\end{itemize}
\end{prop}

An admissible pair $(\sigma,\varepsilon)$ is of trace type if and only if all $(\sigma,\varepsilon)$-sesquilinear forms are trace-valued. By Lemma \ref{trace type pair}, when either $\mathrm{char}(K) \neq 2$ or $\mathrm{char}(K) = 2$ but $\sigma$ acts non-trivially on $Z(K)$, all $(\sigma,\varepsilon)$-sesquilinear forms are trace-valued. 

When $K$ is a field of characteristic $2$ the pair $(\mathrm{id}_K, 1)$ is not of trace type. In this case an $(\mathrm{id}_K, 1)$-sesquilinear form is trace-valued if and only if it is alternating.  

\subsubsection{The polar space $S_f$}\label{Sf}

As in Subsection \ref{Emb-pre}, given a non-zero vector $x\in V$ we denote by $[x]$ the point of $\mathrm{PG}(V)$ represented by the vector $x$ and, for a subspace $X$ of $V$, we set $[X] = \{[x]\}_{x\in X-\{0\}}$. We also write $[x_1, x_2,..., x_k]$ for $[\langle x_1, x_2,..., x_k\rangle]$, for short. 

Given a $(\sigma,\varepsilon)$-sesquilinear form $f:V\times V\rightarrow K$, we say that a point $[x]$ of $\mathrm{PG}(V)$ is {\em isotropic} for $f$ (also $f$-{\em isotropic}) if  the vector $x$ is $f$-isotropic. Similarly, given a subspace $X$ of $V$, the subspace $[X]$ of $\mathrm{PG}(V)$ is {\em totally isotropic} for $f$ ({\em totally $f$-isotropic}) if $X$ is totally $f$-isotropic. We denote by $P_f$ and $L_f$ the set of $f$-isotropic points and totally $f$-isotropic lines of $\mathrm{PG}(V)$ and we put $S_f = (P_f, L_f)$.

Assume that $P_f \neq \emptyset \neq L_f$. Then $S_f$ is a polar space (Buekenhout and Cohen \cite[Chapter 7]{BC}). We call it the polar space {\em associated to} $f$. The singular subspaces of $S_f$ are the totally $f$-isotropic subspaces of $\mathrm{PG}(V)$. The subspace $[\mathrm{Rad}(f)]$ is the radical of $S_f$. So, $S_f$ is non-degenerate if and only if $f$ is non-degenerate.  

The set $P_f$ spans $\mathrm{PG}(V)$ if and only if $f$ is either trivial or trace-valued (Proposition \ref{trace-prop2}, claim (2)).  

Let $e_f:S_f\rightarrow \mathrm{PG}(V)$ be the inclusion mapping of $S_f$ in $\mathrm{PG}(V)$. If $P_f$ spans $\mathrm{PG}(V)$ then $e_f$ is a projective embedding in the sense of Subsection \ref{Emb}.  

\subsubsection{Proportionality of reflexive sesquilinear forms}\label{proportional forms sec}

Let $f:V\times V\rightarrow K$ be a non-trivial $(\sigma,\varepsilon)$-sesquilinear form and let $\kappa\in K-\{0\}$. It is well known (see e.g. Tits \cite[Chapter 8]{Tits}) that $\kappa f$ is a $(\sigma',\varepsilon')$-sesquilinear form where $(\sigma',\varepsilon') = \kappa\cdot(\sigma,\varepsilon)$ (notation as in Subsection \ref{proportional pairs sec}). We say that $f$ and $f'$ are {\em proportional}. 

Clearly, proportional reflexive sesquilinear forms define the same orthogonality relation. A partial converse of this fact also holds, but in order to state it we need one more definition: the {\em non-degenerate rank} of a polar space $S$ is the rank of the quotient of $S$ over its radical (Buekenhout and Cohen \cite[7.5.1]{BC}). The next proposition is implicit in the theory developed in Chapter 9 of Buekenhout and Cohen \cite{BC}.     

\begin{prop}\label{proportional forms prop2}
For $i = 1, 2$, let $(\sigma_i,\varepsilon_i)$ be an admissible pair of $K$ and let $f_i:V\times V\rightarrow K$ be a $(\sigma_i,\varepsilon_i)$-sesquilinear form. Let $S = (P, L)$ be a full subgeometry of $\mathrm{PG}(V)$, satsfying all the following:

\begin{itemize}
\item[$(1)$] The point-set $P$ of $S$ spans $\mathrm{PG}(V)$.
\item[$(2)$] The geometry $S$ is a polar space with non-degenerate rank at least $2$.
\item[$(3)$] The polar space $S$ is a subspace of either of $S_{f_1}$ and $S_{f_2}$.
\end{itemize}
Then the forms $f_1$ and $f_2$ are proportional. 
\end{prop} 

\begin{co}
Given $f_1, f_2:V\times V\rightarrow K$ as in Proposition {\rm {\ref{proportional forms prop2}}}, assume that $S_{f_1} = S_{f_2}$ and the polar space $S := S_{f_1} = S_{f_2}$ has non-degenerate rank at least $2$. Then $f_1$ and $f_2$ are proportional.  
\end{co}

\subsection{Pseudo-quadratic forms}\label{pseudo-quadratic} 

Given a division ring $K$ and an admissible pair $(\sigma,\varepsilon)$ of $K$, let $\overline{K} = K^{(\sigma,\varepsilon)}$, as in (\ref{Kse2}) of Subsection \ref{section 2.1}. The scalar multiplication $\circ$ is defined as in (\ref{circ}) and, for $t\in K$, we write $\bar{t}$ for $t + K_{\sigma,\varepsilon}$, as in Subsection \ref{section 2.1}. 

Let $V$ be a right $K$-vector space. A $(\sigma,\varepsilon)$-{\em quadratic form} on $V$ is a map $q:V\rightarrow\overline{K}$ such that 

\begin{itemize}
\item[(Q1)] ~ $q(x\lambda) = q(x)\circ\lambda$ for any $\lambda\in K$;
\item[(Q2)] there exists a trace-valued $(\sigma,\varepsilon)$-sesquilinear form $f:V\times V\rightarrow K$ such that $q(x+y) = q(x) + q(y) + \overline{f(x,y)}$ for any choice of $x,y\in V$. 
\end{itemize}
We call $f$ a {\em sesquilinearization} of $q$. Note that in the above definition we allow $\overline{K} = \{\bar{0}\}$ (namely $K_{\sigma,\varepsilon} = K$, equivalently $(\sigma,\varepsilon) = (\mathrm{id}_K,-1)$ and $\mathrm{char}(K) \neq 2$), but we warn that when $\overline{K} = \{\bar{0}\}$ both conditions $(Q1)$ and $(Q2)$ are vacuous. In particular, when $\overline{K} = \{\bar{0}\}$ every trace-valued $(\sigma,\varepsilon)$-sesquilinear form satisfies $(Q2)$. On the other hand: 

\begin{lemma}\label{unique sesqui}
Let $\overline{K} \neq \{\bar{0}\}$. Then $q$ admits a unique sesquilinearization.  
\end{lemma}
\pr This lemma is very well known (see Tits \cite[Chapter 8]{Tits}, for instance). Nevertheless, it is worth recalling its proof here, as we shall refer to it later, in Section 3, when discussing generalized pseudo-quadratic forms.

let $f$ and $f'$ be sesquilinearizations of $q$. Then $f(x,y)-f'(x,y) \in K_{\sigma,\varepsilon}$ for any two vectors $x, y\in V$. As 
$f(x\lambda, y\mu) - f'(x\lambda,y\mu) = \lambda^\sigma(f(x,y)-f'(x,y))\mu$
we also have that 
\begin{equation}\label{sesqui-unique}
\lambda^\sigma(f(x,y)-f'(x,y))\mu\in K_{\sigma,\varepsilon} ~~ \mbox{for any} ~ \lambda, \mu\in K.
\end{equation}
If $f(x,y)-f'(x,y) \neq 0$ for some vectors $x, y\in V$ then (\ref{sesqui-unique}) implies that $K_{\sigma,\varepsilon} = K$. However $K_{\sigma,\varepsilon} \subset K$ by assumption. Hence $f(x,y) = f'(x,y)$ for any $x, y\in V$, namely $f = f'$.  \eop 

\bigskip 

In the literature, $(\sigma,\varepsilon)$-quadratic forms are also called {\em pseudo-quadratic forms}, keeping the word {\em quadratic forms} only for $(\mathrm{id}_K, 1)$-quadratic forms.  

We say that a pseudo-quadratic form $q$ is {\em trivial} if $q(x) = \bar{0}$ for any $x\in V$. Clearly, if $\overline{K} = \{\bar{0}\}$ then $q$ is trivial.  

\bigskip

\noindent
{\bf Remark.} In the literature, pseudo-quadratic forms are defined only when $\overline{K}\neq \{\bar{0}\}$. However, in the theory of generalized pseudo-quadratic forms, to be exposed in Section 3, we shall allow forms with trivial codomain. Accordingly, we have allowed $\overline{K} = \{\bar{0}\}$ in our definition of pseudo-quadratic forms.    

\subsubsection{Facilitating forms}\label{facilitating forms quadratic}

Every $(\sigma,\varepsilon)$-quadratic form $q$ admits a so-called {\em facilitating form}, namely a $\sigma$-sesquilinear form $g:V\times V\rightarrow K$ such that 
\begin{equation}\label{facilitating quadratic}
q(x) = \overline{g(x,x)} ~~ \mbox{for any} ~ x\in V.
\end{equation} 
If $\overline{K} = \bar{0}$ every $\sigma$-sesquilinear form is a facilitating form for $q$. Let $\overline{K}\neq \bar{0}$ and let $f$ be the sesquilinearization of $q$. Then all facilitating forms of $q$ are obtained as follows (Tits \cite[Chapter 8]{Tits}). Let $(e_i)_{i\in I}$ be a basis of $V$. Assume that a total ordering $<$ is given on the index set $I$. For every $i\in I$ let $g_i\in K$ be such that $q(e_i) = \bar{g}_i$. For any two vectors $x = \sum_{i\in I}e_i\lambda_i$ and $y = \sum_{i\in I}e_i\mu_i$ of $V$, put 
\begin{equation}\label{facilitating quadratic 1}
g(x, y) := \sum_{i < j}\lambda_i^\sigma f(e_i,e_j)\mu_j + \sum_{i\in I}\lambda_i^\sigma g_i \mu_i.
\end{equation}
(We warn that all sums occurring in (\ref{facilitating quadratic 1}) are well defined, since only finitely many of the scalars $\lambda_i$ and $\mu_i$ are different from $0$.) Then the mapping $g$ defined as in (\ref{facilitating quadratic 1}) is a facilitating form for $q$. Moreover,
\begin{equation}\label{facilitating quadratic 2}
f(x,y) = g(x,y) + g(y,x)^\sigma\varepsilon ~~ \mbox{for any}~ x, y\in V.
\end{equation}
Conversely, given a $\sigma$-sesquilinear form $g:V\times V\rightarrow K$ and an element $\varepsilon\in K$ forming an admissible pair with $\sigma$, let $q:V\rightarrow \overline{K}$ be defined as in (\ref{facilitating quadratic}). Then $q$ is a $(\sigma,\varepsilon)$-quadratic form and the form $f$ defined as in (\ref{facilitating quadratic 2}) is the sesquilinearization of $q$. Note that $f$ is indeed trace-valued, by claim (1) of Proposition \ref{trace-prop2}.

\subsubsection{The polar space $S_q$}

Let $q:V\rightarrow \overline{K}$ be a $(\sigma,\varepsilon)$-quadratic form. We say that a vector $x\in V$ is {\em singular} for $q$ (also $q$-singular) if $q(x) = \bar{0}$. A subspace $X\subset V$ is said to be {\em totally singular} for $q$ (also {\em totally $q$-singular}) if $q(x) =\bar{0}$ for every $x\in X$. 

Clearly, if $q(x) = \bar{0}$ for a vector $x\in V$ then $q(x\lambda) = \bar{0}$ for any $\lambda\in K$. Therefore a point $[x]$ of $\mathrm{PG}(V)$ is totally $q$-singular as a $1$-dimensional subspace of $V$ if and only if $x$ is $q$-singular. If this is the case then we say that the point $[x]$ is {\em singular} for $q$ (also $q$-{\em singular}). A subspace $[X]$ of $\mathrm{PG}(V)$ is said to be {\em totally singular} for $q$ (also {\em totally $q$-singular}) if all of its points are $q$-singular. We denote by $P_q$ and $L_q$ the set of $q$-singular points and totally $q$-singular lines of $\mathrm{PG}(V)$ and we put $S_q = (P_q, L_q)$. 

Note that $P_q$ or $L_q$ could be empty. The opposite situation, where $S_q = \mathrm{PG}(V)$, occurs when $q$ is trivial, as when $\overline{K} = \bar{0}$.  

For the rest of this subsection we assume that $P_q \neq \emptyset \neq L_q$ and $\overline{K}\neq 0$. We denote by $f$ the sesquilinearization of $q$. 

All propositions to be stated in the rest of this subsection are well known. Their proofs can be found in Tits \cite[Chapter 8]{Tits} and Buekenhout and Cohen \cite[Chapter 10]{BC}. However we shall recall those proofs here, since in Section 3 we will need them for reference.    

\begin{prop}\label{singular1 quad}
The point-line geometry $S_q = (P_q,L_q)$ is a subspace of the polar space $S_f$ associated to $f$. Explicitly:
\begin{itemize}
\item[$(1)$] $P_q\subseteq P_f$;
\item[$(2)$] a projective line $[x,y]$ belongs to $L_q$ if and only if $q(x) = q(y) = \bar{0}$ and $f(x,y) = 0$. 
\end{itemize}
\end{prop}
\pr Let $q(x) =\bar{0}$. Then $q(x(\lambda+\mu)) = \bar{0}$ as well, for any choice of scalars $\lambda, \mu \in K$. It follows from $(Q2)$ with $x$ and $y$ replaced by $x\lambda$ and $x\mu$ respectively that $\lambda^\sigma f(x,x)\mu \in K_{\sigma,\varepsilon}$ for any choice of $\lambda$ and $\mu$. If $f(x,x)\neq 0$, this forces $K_{\sigma,\varepsilon} = K$, contradicting the assumption that $\overline{K}\neq \bar{0}$. Therefore $f(x,x) = 0$. Claim (1) is proved. 

Turning to claim (2), let $[x,y]\in L_q$. Then $q(x\lambda + y\mu) = \bar{0}$ for any choice of $\lambda, \mu \in K$. According to $(Q2)$, this forces $\lambda^\sigma f(x,y)\mu \in K_{\sigma,\varepsilon}$ for all $\lambda, \mu\in K$. Hence $f(x,y) = 0$, since $K_{\sigma, \varepsilon}\subset K$. The `only if' part of (2) is proved. The `if' part is trivial.  \eop 

\bigskip

The next two corollaries immediately follow from Proposition \ref{singular1 quad}. 

\begin{co}\label{singular2 quad} 
A subspace $[x_1, x_2,..., x_k]$ of $\mathrm{PG}(V)$ is totally $q$-singular if and only if it is totally isotropic for $f$ and
$q(x_1) = q(x_2) = ... = q(x_k) = \bar{0}$.
\end{co}

\begin{co}\label{singular 3 quad}
The point-line geometry $S_q$ is a polar space. Its singular subspaces are the totally $q$-singular subspaces of $\mathrm{PG}(V)$. The set $P_q\cap[\mathrm{Rad}(f)]$ is the radical of $S_q$. 
\end{co} 

The radical $P_q\cap[\mathrm{Rad}(f)]$ of $S_q$ is a subspace of $[\mathrm{Rad}(f)]$. In other words, the 
$q$-singular vectors of $\mathrm{Rad}(f)$ form a subspace of $\mathrm{Rad}(f)$. We call this subspace the {\em radical} of $q$ and we denote it by the symbol $\mathrm{Rad}(q)$.  Following Buekenhout and Cohen \cite[Chapter 10]{BC} we call $\mathrm{Rad}(f)$ the {\em defect} of $q$ (but we warn that this word is used with a different meaning in Tits \cite{Tits}). 

The form $q$ is said to be {\em singular} (also {\em degenerate})  if $\mathrm{Rad}(q) \neq \{0\}$. 

\begin{prop}\label{Sq in Sf generation}
If $P_q\not\subseteq [\mathrm{Rad}(f)]$ then $P_q$ spans $\mathrm{PG}(V)$.
\end{prop}
\pr Suppose that $P_q\not\subseteq[\mathrm{Rad}(f)]$. Then there exists a $q$-singular point $[a]\not\in [\mathrm{Rad}(f)]$. As $a\not\in \mathrm{Rad}(f)$, the space $a^\perp$ is a hyperplane of $V$. Let $l = [a,b]$ be a projective line of $\mathrm{PG}(V)$ through $[a]$ not contained in $[a^\perp]$. Then $f(a,b) \neq 0$. Moreover, 
\begin{equation}\label{span eq1}
q(a\lambda + b) = q(a)\circ\lambda + q(b) + \overline{\lambda^\sigma f(a,b)} = q(b) +  \overline{\lambda^\sigma f(a,b)}
\end{equation}
by $(Q2)$ and since $q(a) =\bar{0}$. As $f(a,b)\neq 0$, there exists a scalar $\lambda\in K$ such that $q(b) + \overline{\lambda^\sigma f(a,b)} = \bar{0}$. Then $q(a\lambda+b) = \bar{0}$ by (\ref{span eq1}). So, the vector $b_l := a\lambda+b$ is $q$-singular and $[b_l]\neq [a]$. 

Let $\Lambda_a$ be the set of lines of $\mathrm{PG}(V)$ that contain $[a]$ but are not contained in $[a^\perp]$. By the previous paragraph, every line $l\in \Lambda_a$ contains a $q$-singular point $[b_l] \neq [a]$. Let $\Pi_a := \{[b_l]\}_{l\in \Lambda_a}$. Then $\Pi_a$ is contained in $P_q$ and spans $\mathrm{PG}(V)$. Hence $\langle P_q\rangle = \mathrm{PG}(V)$. \eop 

\bigskip  

If $P_q$ spans $\mathrm{PG}(V)$ then the inclusion mapping $e_q:S_q\rightarrow \mathrm{PG}(V)$ is a projective embedding in the sense of Subsection \ref{Emb}. 

We know that $S_q$ is a subspace of $S_f$ (Proposition \ref{singular1 quad}), but it could be a proper subspace of $S_f$, namely vectors $x\in V$ might exist such that $f(x,x) = 0$ but $q(x)\neq \bar{0}$. Nevertheless, the following holds.  

\begin{lemma}\label{from iso to q}
For $x\in V$, if $f(x,x) = 0$ then $q(x) \in \overline{K}^\circ$.
\end{lemma}
\pr Recall that $\overline{K}^\circ = K^{\sigma,\varepsilon}/K_{\sigma,\varepsilon}$ (Lemma \ref{lemma1}). Let $f(x,x) = 0$. Then 
\[\begin{array}{rcl}
q(x)\circ(\lambda+\mu) & =&  q(x(\lambda+\mu))  ~ = ~  q(x\lambda) + q(x\mu) + \overline{\lambda^\sigma f(x,x)\mu} ~ =\\
{} & = &  q(x\lambda) + q(x\mu)  ~ = ~ q(x)\circ\lambda +q(x)\circ\mu
\end{array}\]
for any choice of $\lambda, \mu \in K$. Let $t\in K$ be such that $q(x) = \bar{t}$. By the above, we have $(\lambda+\mu)^\sigma t(\lambda+\mu) \equiv \lambda^\sigma t\lambda + \mu^\sigma t\mu ~~ (\mathrm{mod} ~ K_{\sigma,\varepsilon})$. Hence
\begin{equation}\label{f(x,x)=0;1}
\lambda^\sigma t\mu + \mu^\sigma t\lambda \in K_{\sigma,\varepsilon}.
\end{equation} 
Recalling that $\lambda^\sigma t\mu - (\lambda^\sigma t\mu)^\sigma\varepsilon \in K_{\sigma,\varepsilon}$ and $(\lambda^\sigma t\mu)^\sigma\varepsilon = \mu^\sigma t^\sigma\varepsilon \lambda$, from (\ref{f(x,x)=0;1}) we obtain that $\mu^\sigma t\lambda + \mu^\sigma t^\sigma\varepsilon\lambda \in K_{\sigma,\varepsilon}$, namely
\begin{equation}\label{f(x,x)=0;2}
\mu^\sigma(t + t^\sigma\varepsilon)\lambda  \in K_{\sigma,\varepsilon}.
\end{equation} 
As $K_{\sigma,\varepsilon} \neq K$ by assumption and (\ref{f(x,x)=0;2}) holds for any choice of $\lambda, \mu\in K$, we obtain that $t + t^\sigma\varepsilon = 0$, namely $t \in K^{\sigma,\epsilon}$. Hence $\bar{t}\in K^{\sigma,\varepsilon}/K_{\sigma,\varepsilon} = \overline{K}^\circ$.  \eop 

\begin{prop}\label{Sq trace type}
Let $(\sigma,\varepsilon)$ be of trace type. Then $S_q = S_f$.
\end{prop}
\pr Let $(\sigma,\varepsilon)$ be of trace type. Then $\overline{K}^\circ =\bar{0}$, by claim (1) of Corollary \ref{cor1}. Lemma \ref{from iso to q} now implies that all $f$-isotropic vectors are $q$-singular, namely $P_f \subseteq P_q$. Therefore $S_q = S_f$, since $S_q$ is a subspace of $S_f$.  \eop  

\subsubsection{Proportionality of pseudo-quadratic forms}\label{proportional quadric sec}

In this subsection we adopt the notation of Subsection \ref{proportional pairs sec}, thus denoting the group $\overline{K} = K/K_{\sigma,\varepsilon}$ by the symbol $\overline{K}^{\sigma,\varepsilon}$. 

Assuming that $K_{\sigma,\varepsilon} \neq K$, let $q:V\rightarrow\overline{K}^{\sigma,\varepsilon}$ be a non-trivial $(\sigma,\varepsilon)$-quadratic form and let $f$ be its sesquilinearization. Given a
scalar $\kappa\in K-\{0\}$, let $(\sigma',\varepsilon') := \kappa\cdot(\sigma,\varepsilon)$. Let $\kappa q:V \rightarrow \overline{K}^{\sigma',\varepsilon'}$ map every $x\in V$ onto $\kappa q(x)\in \overline{K}^{\sigma',\varepsilon'}$ (well defined by Lemma \ref{proportional pairs lemma}). Then $\kappa q$ is a $(\sigma',\varepsilon')$-quadratic form and $\kappa f$ is the sesquilinearization of $\kappa q$ (Tits \cite[Chapter 8]{Tits}). Clearly, $S_{q'} = S_q$. We say that $q$ and $q'$ are {\em proportional}.  

\begin{prop}\label{proportional quadric prop}
For $i = 1, 2$, let $q_i:V\rightarrow \overline{K}^{\sigma_i,\varepsilon_i}$ be a non-trivial $(\sigma_i,\varepsilon_i)$-quadratic form such that $S_{q_i}$ has non-degenerate rank at least $2$. Suppose that $S_{q_1} = S_{q_2}$. Then $q_1$ and $q_2$ are proportional. 
\end{prop}
\pr This proposition is well known (see e.g. Tits \cite[Chapter 8]{Tits}). Nevertheless we give a sketch of the proof here,  since in the Section 3 we will need it for reference. 

Let $f_1$ and $f_2$ be the sesquilinearizations of $q_1$ and $q_2$. By Proposition \ref{Sq in Sf generation}, for $i = 1, 2$ the set $P_{q_i}$ spans $\mathrm{PG}(V)$. Moreover $S_{q_i}$ is a subspace of $S_{f_i}$. By assumption, the polar space $S_{q_i}$ has non-degenerate rank at least $2$. Hence the equality $S_{q_1} = S_{q_2}$ forces $f_1$ and $f_2$ to be proportional, by Proposition \ref{proportional forms prop2}. It follows that $q_1$ and $q_2$ admit proportional facilitating forms (see definition (\ref{facilitating quadratic 1}), with a basis of singular vectors). Hence they are proportional.  \eop     

\section{Generalized pseudo-quadratic forms} 

In this section we propose a generalization of pseudo-quadratic forms and we show that all what we have said on the latters in the previous section remains valid in this more general context. 

\subsection{Definition and basic properties} 

Given a division ring $K$ and an admissible pair $(\sigma,\varepsilon)$ of $K$, let $\overline{R}$ be a $\circ$-closed subgroup of $\overline{K}$ (see Subsection \ref{scalar-mult-sec}). We denote by $R$ the pre-image of $\overline{R}$ by the projection $t \mapsto \bar{t} = t+ K_{\sigma,\varepsilon}$ of $K$ onto $\overline{K} = K/K_{\sigma,\varepsilon}$:
\begin{equation}\label{pre}
R := \{t ~|~ \bar{t}\in\overline{R}\}.
\end{equation}
We recall that a scalar multiplication is induced by $\circ$ on the factor group $\overline{K}/\overline{R}$, as explained in (\ref{circ-quot}). Clearly $\overline{R}$ is the null element of $\overline{K}/\overline{R}$. When $\overline{R}$ is given this role, we denote it by the symbol $0_{\overline{R}}$. 

Given a $K$-vector space $V$, a {\em generalized} $(\sigma,\varepsilon)$-{\em quadratic form} (also {\em generalized pseudo-quadratic form}) is a map $q:V\rightarrow\overline{K}/\overline{R}$ such that 

\begin{itemize}
\item[(Q'1)] ~ $q(x\lambda) = q(x)\circ\lambda$ for any $\lambda\in K$;
\item[(Q'2)] there exists a trace-valued $(\sigma,\varepsilon)$-sesquilinear form $f:V\times V\rightarrow K$ such that $q(x+y) = q(x) + q(y) + (\overline{f(x,y)}+\overline{R})$ for any choice of $x,y\in V$.
\end{itemize}
We call $\overline{R}$ the {\em co-defect} of $q$. With this terminology, a pseudo-quadratic form is just a generalized pseudo-quadratic form with trivial co-defect. 

\bigskip

\noindent
{\bf Remark.} A motivation for the choice of the word {\em co-defect} will be given in Subsection \ref{dominant covers}, where we will show that the co-defect  $\overline{R}$ of $q$ is involved in the defect of a suitable pseudo-quadratic form, called the dominant cover of $q$.

\bigskip 

A sesquilinear form $f$ as in $(Q'2)$ is called a {\em sesquilinearization} of $q$. 

\begin{lemma}\label{form}
Let $q:V\rightarrow \overline{K}/\overline{R}$ be a generalized pseudo-quadratic form. 

\begin{itemize}
\item[$(1)$] If $\overline{R}\neq \overline{K}$ then $q$ admits exactly one sesquilinearization. 
\item[$(2)$] Let $\overline{R} = \overline{K}$. Then every trace-valued $(\sigma,\varepsilon)$-sesquilinar form on $V$ is a sesquilinearization of $q$.
\end{itemize}
\end{lemma}
\pr Claim $(2)$ is obvious. Claim $(1)$ can be proved by the same argument used to prove Lemma \ref{unique sesqui}, but for replacing $K_{\sigma,\varepsilon}$ with the group $R$ defined in (\ref{pre}).   \eop  

\bigskip

Every generalized $(\sigma,\varepsilon)$-quadratic form also admits a {\em facilitating form}, namely a $\sigma$-sesquilinear form $g:V\times V\rightarrow K$ such that 
\begin{equation}\label{facilitating gen quadratic}
q(x) = \overline{g(x,x)} + \overline{R} ~~ \mbox{for any} ~ x\in V.
\end{equation} 
If $\overline{R} = \overline{K}$ then every $\sigma$-sesquilinear form is a facilitating form for $q$. Let $\overline{R}\neq\overline{K}$ and let $f$ be the sesquilinearization of $q$. It is straightforward to prove that all facilitating forms of $q$ are obtained as follows. Let $(e_i)_{i\in I}$ be a basis of $V$ and $<$ a total ordering of $I$. For every $i\in I$ let $g_i\in K$ be such that $q(e_i) = \bar{g}_i + \overline{R}$. For $x,y\in V$ let $g(x,y)$ be defined as in (\ref{facilitating quadratic 1}). Then $g$ is a facilitating form for $q$. Moreover $f(x,y) = g(x,y) + g(y,x)^\sigma\varepsilon$, as in (\ref{facilitating quadratic 2}). 

Conversely, given a $\sigma$-sesquilinear form $g:V\times V\rightarrow K$ and an element $\varepsilon\in K$ forming an admissible pair with $\sigma$, let $q:V\rightarrow \overline{K}$ be defined as in (\ref{facilitating gen quadratic}). Then $q$ is a generalized $(\sigma,\varepsilon)$-quadratic form and the form $f$ defined as in (\ref{facilitating quadratic 2}) is the sesquilinearization of $q$.   

\begin{theo}\label{forms-theo}
Let $\overline{R} \neq \overline{K}$. Let $q:V\rightarrow \overline{K}/\overline{R}$ be a generalized $(\sigma,\varepsilon)$-quadratic form, let $f$ be its sesquilinearization and $R$ as in {\rm (\ref{pre})}. Then all the following hold:
\begin{itemize}
\item[(1)] We have $\overline{R}\subseteq \overline{K}^\circ$. In other words, $\overline{R}$ is a vector subspace of $\overline{K}^\circ$. 
\item[(2)] For every vector $x\in V$, if $q(x) = 0_{\overline{R}}$ then $f(x,x) = 0$.
\item[(3)] Let $x\in V$ be such that $f(x,x) = 0$. Then $q(x)\in \overline{K}^\circ/\overline{R}$ (well defined in view of claim $(1)$). 
\end{itemize}
\end{theo}
\pr In view of (Q'1) and (Q'2), we have 
\[q(x)\circ(\lambda+\mu) + \overline{R} = q(x(\lambda+\mu)) = q(x)\circ\lambda + q(x)\circ\mu + \lambda^\sigma f(x,x)\mu\] 
for any choice of $\lambda,\mu\in K$. Therefore, given $t\in K$ such that $\bar{t} + \overline{R} = q(x)$, we have
$\lambda^\sigma t\mu + \mu^\sigma t\lambda - \lambda^\sigma f(x,x)\mu \in R$. As $K_{\sigma,\varepsilon}\subseteq R$ and $\mu^\sigma t\lambda - \lambda^\sigma t^\sigma\varepsilon \mu = \mu^\sigma t\lambda - (\mu^\sigma t\lambda)^\sigma\varepsilon \in K_{\sigma,\varepsilon}$ 
we obtain that $\lambda^\sigma t\mu + \lambda^\sigma t^\sigma\varepsilon\mu- \lambda^\sigma f(x,x)\mu \in R$, namely
\begin{equation}\label{R0}
\lambda^\sigma(t+t^\sigma\varepsilon - f(x,x))\mu \in R ~~\mbox{for any choice of} ~ \lambda, \mu \in K.
\end{equation}
As $R \subset K$ by assumption, (\ref{R0}) forces 
\begin{equation}\label{R1}
t+t^\sigma\varepsilon = f(x,x).
\end{equation}
However we can replace $t$ with $t+r$ in (\ref{R1}), for any $r\in R$. By comparing the new equation thus obtained with (\ref{R1}) we obtain that $r+r^\sigma\varepsilon = 0$ for any $r\in R$, namely $R\subseteq K^{\sigma,\varepsilon}$. Equivalently, $\overline{R}\subseteq K^{\sigma,\varepsilon}/K_{\sigma,\varepsilon} = \overline{K}^\circ$, as claimed in (1). As $\overline{R}$ is $\circ$-closed by assumption, $\overline{R}$ is a vector subspace of the $K$-vector space $\overline{K}^\circ$.  

Claims (2) and (3) can be proved in the same way as claim (1) of Proposition \ref{singular1 quad} and Lemma \ref{from iso to q}, but for replacing $K_{\sigma,\varepsilon}$ with $R$ in those proofs.   \eop 

\bigskip

Note that $f(x,x)\in K_{\sigma,-\varepsilon}$ for any $x\in V$ because $f$ is trace-valued. If $\mathrm{char}(K) = 2$ then $\varepsilon = -\varepsilon$. In this case $f(x,x) \in K_{\sigma,\varepsilon} \subseteq R$ for any $x\in V$.  

\begin{co}\label{forms-cor}
Let $(\sigma,\varepsilon)$ be of trace type and $\overline{R}\neq \overline{K}$. Then $\overline{R} =\{\bar{0}\}$, whence $q$ is pseudo-quadratic. 
\end{co}
\pr By claim (1)  of Corollary \ref{cor1}, the pairs $(\sigma,\varepsilon)$ is of trace type if and only if $\overline{K}^\circ = \{\bar{0}\}$. Moreover, by claim (1) of Theorem \ref{forms-theo}, either $\overline{R} = \overline{K}$ or $\overline{R}\subseteq\overline{K}^\circ$. Therefore, if $\overline{R} \subset \overline{K}$ and $\overline{K}^\circ = \{\bar{0}\}$ then $\overline{R} = \{\bar{0}\}$.  \eop 

\bigskip

A generalized pseudo-quadratic form $q:V\rightarrow \overline{K}/\overline{R}$ is said to be {\em trivial} if $q(x) = 0_{\overline{R}}$ for every $x\in V$.  

\begin{prop}\label{triviality}
The form $q$ is trivial if and only if one of the following holds:
\begin{itemize}
\item[$(1)$] $\overline{R} = \overline{K}$.
\item[$(2)$] We have $\overline{R} \neq \overline{K}$ but the sesquilinearization of $q$ is trivial and there exists a basis $(e_i)_{i\in I}$ of $V$ such that $q(e_i) = 0_{\overline{R}}$ for every $i\in I$.
\end{itemize}
\end{prop}
\pr Clearly, if $\overline{R} = \overline{K}$ then $q$ is trivial. Assume that $\overline{R}\subset\overline{K}$. Then $q$ admits a unique sesquilinearization $f$, by Lemma \ref{form}. Suppose that nevertheless $q$ is trivial. Then $f(x,y)\in R$ for any $x, y\in V$. Accordingly, 
\begin{equation}\label{trivial}
\lambda^\sigma f(x,y)\mu \in R ~~ \mbox{for any choice of} ~ \lambda, \mu\in K ~\mbox{and} ~ x, y\in  V. 
\end{equation}
If $f(x,y)\neq 0$ for a pair $(x,y)$, then (\ref{trivial}) forces $R = K$, contrary to the assumptions made on $\overline{R}$. It follows that $f(x,y)$ is the trivial form.  

Conversely, let $f$ be trivial and $q(e_i) = 0_{\overline{R}}$ for every $i \in I$. Then the form $g$ defined as in (\ref{facilitating quadratic 2}) but with $g_i = 0$ for every $i\in I$, is trivial. However $g$ is a facilitating form of $q$. Hence $q$ is trivial as well.  \eop 

\subsection{The polar space $S_q$}\label{Sq gen quadr} 

For the rest of this section we assume that $q$ is non-trivial. In particular, $\overline{R}\neq \overline{K}$. As above, $f$ stands for the sesquilinearization of $q$. The symbol $R$ is given the meaning stated in (\ref{pre}). 

As in the case of pseudo-quadratic forms, we say that a vector $x\in V$ is {\em singular} for $q$ (also $q$-{\em singular}) if $q(x) = 0_{\overline{R}}$. A subspace $X$ of $V$ is said to be {\em totally singular} for $q$ (also {\em totally $q$-singular}) if $q(x) = 0_{\overline{R}}$ for every $x\in X$. 

Clearly, if $q(x) = 0_{\overline{R}}$ for a vector $x\in V$ then $q(x\lambda) = 0_{\overline{R}}$ for any $\lambda\in K$. We say that a point $[x]$ of $\mathrm{PG}$ is $q$-{\em singular} (also $q$-{\em singular}) if $x$ is $q$-singular. A subspace of $\mathrm{PG}(V)$ is said to be {\em totally singular} for $q$ ({\em totally $q$-singular}) if all of its points are $q$-singular. 

Let $P_q$ be the set of $q$-singular points of $\mathrm{PG}(V)$. By claim (2) of Theorem \ref{forms-theo}, if a point of $\mathrm{PG}(V)$ is $q$-singular then it is $f$-isotropic. In short, $P_q\subseteq P_f$.

\begin{prop}\label{singular1}
A line $[x,y]$ of $\mathrm{PG}(V)$ is totally $q$-singular if and only if $q(x) = q(y) = 0_{\overline{R}}$ and $f(x,y) = 0$. 
\end{prop}
\pr This statement can be proved in the same way as claim (2) of Proposition \ref{singular1 quad}, but for replacing $K_{\sigma,\varepsilon}$ with $R$ in that proof.   \eop  

\begin{co}\label{singular2} 
A subspace $[x_1, x_2,..., x_k]$ of $\mathrm{PG}(V)$ is totally $q$-singular if and only if it is totally isotropic for $f$ and
$q(x_1) = q(x_2) = ... = q(x_k) = 0_{\overline{R}}$.
\end{co}
\pr This immediately follows from Proposition \ref{singular1}.  \eop

\begin{co}\label{singular2-bis}
Let $(\sigma,\varepsilon)$ be of trace type. Then a subspace of $\mathrm{PG}(V)$ is totally $q$-singular if and only it is totally $f$-isotropic. 
\end{co}
\pr By Corollary \ref{forms-cor}, when $(\sigma,\varepsilon)$ is of trace type the form $q$ is pseudo-quadratic. The conclusion follows from Proposition \ref{Sq trace type}. \eop 

\bigskip

Assuming that $P_q\neq \emptyset$, let $L_q$ be the set of totally $q$-singular lines of $\mathrm{PG}(V)$ and put $S_q = (P_q,L_q)$. In view of Proposition \ref{singular2}, the point-line geometry $S_q$ is a subspace of the polar space $S_f = (P_f,L_f)$ associated to $f$. It readily follows that $S_q$ is itself a polar space. Its radical is a (possibly empty) subspace of $[\mathrm{Rad}(f)]$, equal to $P_q\cap [\mathrm{Rad}(f)]$. Moreover, if $(\sigma,\varepsilon)$ is of trace type then $S_q = S_f$, by Corollary \ref{singular2-bis}. 

We call $S_q$ the polar space {\em associated to} $q$. The $q$-singular vectors of $\mathrm{Rad}(f)$ form a subspace of $\mathrm{Rad}(f)$, henceforth called the {\em radical} of $q$ and denoted by the symbol $\mathrm{Rad}(q)$. We say that $q$ is {\em singular} (also {\em degenerate}) if $\mathrm{Rad}(q) \neq \{0\}$, namely $S_q$ is degenerate. We call $\mathrm{Rad}(f)$ the {\em defect} of $q$, let $q$ be singular or not. 

Let $q_{|\mathrm{Rad}(f)}$ be the mapping induced by $q$ on $\mathrm{Rad}(f)$. Clearly $q_{|\mathrm{Rad}(f)}$ is additive. By this fact and claim (3) of Theorem \ref{forms-theo} we get the following:

\begin{prop}\label{hom}
The mapping $q_{|\mathrm{Rad}(f)}$ is a homomorphism of $K$-vector spaces from $\mathrm{Rad}(f)$ to $\overline{K}^\circ/\overline{R}$ and $\mathrm{Rad}(q)$ is the kernel of this homomorphism.
\end{prop}

Consequently, the quotient space $\mathrm{Rad}(f)/\mathrm{Rad}(q)$ is isomorphic to the image $\mathrm{Im}(q_{|\mathrm{Rad}(f)})$ of $q_{|\mathrm{Rad}(f)}$ and the latter is a vector subspace of $\overline{K}^\circ/\overline{R}$. 

\bigskip

\noindent
{\bf Remark.} A result similar to Proposition \ref{hom} holds with $\mathrm{Rad}(f)$ replaced by any totally $f$-isotropic subspace $X$ of $V$ and $\mathrm{Rad}(q)$ replaced by the set of $q$-singular vectors of $X$.

\begin{prop}\label{span}
Either $P_q$ is totally $q$-singular or it spans $\mathrm{PG}(V)$.
\end{prop}
\pr The proof given for Proposition \ref{Sq in Sf generation} works for this statement as well, but for replacing $K_{\sigma,\varepsilon}$ with $R$ in that proof.   \eop  

\bigskip

When $P_q$ spans $\mathrm{PG}(V)$ the inclusion mapping $e_q:S_q\rightarrow \mathrm{PG}(V)$ is an embedding as defined in Subsection \ref{Emb}. 

\subsection{A facilitating form}\label{facilitating form}

We keep the hypotheses and the notation of the previous subsection. In particular, $\overline{R}\neq \overline{K}$, $f$ is the sesquilinearization of $q$ and $P_q$ is the set of $q$-singular points of $\mathrm{PG}(V)$. We also assume that $P_q$ spans $\mathrm{PG}(V)$. Hence $V$ admits a basis formed by $q$-singular vectors. We call such a basis a $q$-{\em singular basis}. 

Let $E = (e_i)_{i\in I}$ be a $q$-singular basis of $V$. Given a total ordering $<$ on the set $I$ of indices, let $g_E:V\times V\rightarrow K$ be the $\sigma$-sequilinear form defined as follows:
\begin{equation}\label{facilitating-E1}
g_E(\sum_ie_i\lambda_i, \sum_je_j\mu_j) := \sum_{i<j}\lambda^\sigma f(e_i,e_j)\mu_j.
\end{equation}
Since $q(e_i) = 0_{\bar{R}}$ for every $i\in I$, the form $g_E$ is a facilitating form for $q$, namely  
\[q(x) = \overline{g_E(x,x)} + \overline{R} = \sum_{i<j}\overline{\lambda_i^\sigma f(e_i,e_j)\lambda_j} + \overline{R} \]
for every vector $x = \sum_{i\in I}e_i\lambda_i$ of $V$. Clearly, the coset $\overline{g_E(x,x)}+ \overline{R}$ does not depend on the choice of the $q$-singular basis $E$ but the scalar $g_E(x,x)$ obviously depends on that choice. The value $\overline{g_E(x,x)}$ also depends on it, to some extent. In order to make this remark less vague, we need a few additional definitions. 

Let $E = (e_i)_{i\in I}$ and $E' = (e'_i)_{i\in I}$ be two ordered $q$-singular bases of $V$. Let $\overline{R}_{E,E'}$ be the $\circ$-closed subgroup of $\overline{K}$ spanned by the family $\{\overline{g_{E'}(e_i,e_i)}\}_{i\in I}$ and let $\delta_{E,E'}:V\in \overline{K}$ be the mapping defined as follows:
\[\delta_{E,E'}(x) := \overline{g_E(x,x)}-\overline{g_{E'}(x,x)}.\]
Clearly, $\delta_{E,E'}(x) +\overline{R} = q(x)-q(x) = 0_{\overline{R}}$. Therefore $\delta_{E,E'}(V)\subseteq\overline{R}$. Recall that $\overline{R}$ is a vector subspace of $\overline{K}^\circ$, as we know from Theorem \ref{forms-theo}, (1).  

\begin{lemma}\label{irrelevant}
The group $\overline{R}_{E,E'}$, equipped with the scalar multiplication $\circ$, is a vector subspace of $\overline{R}$ and $\delta_{E,E'}$ is a surjective linear map from $V$ to $\overline{R}_{E,E'}$. Moreover $\delta_{E',E} = -\delta_{E.E'}$ and $\overline{R}_{E,E'} = \overline{R}_{E',E}$.  
\end{lemma}
\pr For $x \in V$ let $x  =  \sum_ie_i\lambda_i = \sum_ie'_i\lambda'_i$. Then
\begin{equation}\label{sommona1}
\left.\begin{array}{rcl}
g_E(x,x) & = & \sum_{i<j}\lambda_i^\sigma f(e_i,e_j)\lambda_j,\\
g_{E'}(x,x) & = & \sum_{i<j}(\lambda'_i)^\sigma f(e'_i,e'_j)\lambda'_j.
\end{array}\right\}
\end{equation}
Moreover, there exists scalars $\alpha_{ij}$ $(i,j\in I)$ such that
\begin{equation}\label{sommona2}
e_k ~ = ~ \sum_ie'_i\alpha_{ik} ~~ \mbox{for all}~ k\in I.
\end{equation}
Hence 
\begin{equation}\label{sommona3}
\lambda'_k ~ = ~ \sum_i\alpha_{ki}\lambda_i ~~ \mbox{for all}~ k\in I.
\end{equation} 
Substituting (\ref{sommona2}) in the first equality of (\ref{sommona1}) and (\ref{sommona3}) in the second one we get
\begin{equation}\label{sommona4}
\left.\begin{array}{rcl}
g_E(x,x) & = & \sum_{i<j}\sum_{k,h}\lambda_i^\sigma\alpha_{k,i}^\sigma f(e'_k,e'_h)\alpha_{h,j}\lambda_j,\\
g_{E'}(x,x) & = & \sum_{i<j}\sum_{k,h}\lambda_k^\sigma\alpha_{i,k}^\sigma f(e'_i,e'_j)\alpha_{j,h}\lambda_h.
\end{array}\right\}
\end{equation}
By changing indices in the second equation of (\ref{sommona4}), we can rewrite the two equations of (\ref{sommona4}) as follows:
\begin{equation}\label{sommona5}
\left.\begin{array}{rcl}
g_E(x,x) & = & \sum_{i,j,k,h; ~ i<j}\lambda_i^\sigma\alpha_{k,i}^\sigma f(e'_k,e'_h)\alpha_{h,j}\lambda_j,\\
g_{E'}(x,x) & = & \sum_{i,j,k,h; ~ k<h}\lambda_i^\sigma\alpha_{k,i}^\sigma f(e'_k,e'_h)\alpha_{h,j}\lambda_j.
\end{array}\right\}
\end{equation}
Recalling that $f(e'_h,e'_k)  =  f(e'_k,e'_h)^\sigma\varepsilon$, that
\[\begin{array}{ll}
\lambda_i^\sigma\alpha_{k,i}^\sigma f(e'_k,e'_h)\alpha_{h,j}\lambda_j -
\lambda_j^\sigma\alpha_{h,j}^\sigma f(e'_k,e'_h)^\sigma\varepsilon\alpha_{k,i}\lambda_i = & \\
= (\lambda_i^\sigma\alpha_{k,i}^\sigma f(e'_k,e'_h)\alpha_{h,j}\lambda_j) -
(\lambda_i^\sigma\alpha_{k,i}^\sigma f(e'_k,e'_h)\alpha_{h,j}\lambda_j)^\sigma\varepsilon  & 
\in K_{\sigma,\varepsilon},
\end{array}\]
and $f(e'_k,e'_k) = 0$ (by $(2)$ of Theorem \ref{forms-theo} and since $q(e'_k) = 0_{\overline{R}}$ by assumption), we can rewrite the two equalities of (\ref{sommona5}) as follows: 
\[\begin{array}{rcl}
g_E(x,x) & = &  \sum_{i<j, k<h}\lambda_i^\sigma\alpha_{k,i}^\sigma(f(e'_k,e'_h)+f(e'_k,e'_h)^\sigma\varepsilon)\alpha_{h,j}\lambda_j, \\
& & \\
g_{E'}(x,x) + K_{\sigma,\varepsilon} & = & \sum_{i<j, k<h}\lambda_i^\sigma\alpha_{k,i}^\sigma(f(e'_k,e'_h)+f(e'_k,e'_h)^\sigma\varepsilon)\alpha_{h,j}\lambda_j +\\
 & & +  \sum_{k,h,i; k<h}\lambda_i^\sigma\alpha^\sigma_{k,i}f(e'_k,e'_h)\alpha_{h,i}\lambda_i. \\
\end{array}\]
Consequently,
\begin{equation}\label{sommona6}
\overline{g_E(x,x)} -\overline{ g_{E'}(x,x)} ~ = ~ - \sum_{k,h,i;k<h}\overline{\alpha_{k,i}^\sigma f(e'_k,e'_h)\alpha_{h,i}}\circ\lambda_i.
\end{equation}
However $\sum_{k<h}\alpha_{k,i}^\sigma f(e'_k,e'_k)\alpha_{h,i} = g_{E'}(\sum_ke'_k\alpha_{k,i}, \sum_ke'_k\alpha_{k,i}) = g_{E'}(e_i,e_i)$
by definition of $g_{E'}$ and (\ref{sommona2}). Substituing in (\ref{sommona6}) we obtain:
\begin{equation}\label{sommona7}
\overline{g_E(x,x)} -\overline{ g_{E'}(x,x)} ~ = ~ - \sum_i\overline{g_{E'}(e_i,e_i)}\circ\lambda_i.
\end{equation}
According to (\ref{sommona7}), we have $\overline{R}_{E,E'} = \delta_{E,E'}(V)$ ($\subseteq \overline{R}$, as previously remarked). Therefore $\overline{R}_{E,E'}$ is a vector subspace of $\overline{R}$. Equation (\ref{sommona7}) also shows that $\delta_{E,E'}$ is a linear mapping from $V$ to $\overline{R}_{E,E'}$. Clearly, $\delta_{E',E} = -\delta_{E,E'}$. Whence $\overline{R}_{E,E'} = \overline{R}_{E',E}$.  \eop 

\bigskip

We call $\delta_{E,E'}$ and $\overline{R}_{E,E'}$ the {\em difference-map} and the {\em difference-space} relative to the pair $(E,E')$ of $q$-singular bases. 

\bigskip

\noindent
{\bf Remark.} Only $q$-singular bases are considered in Lemma \ref{irrelevant}, but the statement of Lemma \ref{irrelevant} holds for any pair of bases formed by $f$-isotropic vectors, except that in this more general setting no closed subgroup $\overline{R}$ is given in advance. Instead of $\overline{R}$ we must consider the closed subgroups $\overline{R}_E$ and $\overline{R}_{E'}$ of $\overline{K}$ generated by the sets $\{\overline{g_E(x,x)}\}_{[x]\in P_f}$ and $\{\overline{g_{E'}(x,x)}\}_{[x]\in P_f}$ respectively. The proof of Lemma \ref{irrelevant} shows that $\delta_{E,E'}(V) = \overline{R}_{E,E'}\subseteq \overline{R}_{E'}$, whence $\overline{R}_E\subseteq\overline{R}_{E'}$. By symmetry, $\overline{R}_{E}\supseteq \overline{R}_{E'}$. Finally $\overline{R}_{E} = \overline{R}_{E'}$.  

For every $x\in V$, put $\gamma_E(x) := \overline{g_E(x,x)}$ and $\gamma_{E'}(x) := \overline{g_{E'}(x,x)}$. Then both $\gamma_E$ and $\gamma_{E'}$ are pseudo-quadratic forms. By Lemma \ref{from iso to q}, the group $\overline{R}_E = \overline{R}_{E'}$ is a vector subspace of $\overline{K}^\circ$. 

\subsection{Isomorphism and weak isomorphism}\label{weak iso}

Given two generalized $(\sigma,\varepsilon)$-quadratic forms $q:V\rightarrow \overline{K}/\overline{R}$ and $q':V'\rightarrow \overline{K}/\overline{R}$ with the same co-defect $\overline{R}$, we say that $q$ and $q'$ are {\em isomorphic} if there exists a bijective linear mapping $\alpha:V\rightarrow V'$ such that $q'(\alpha(x)) = q(x)$ for every $x\in V$. 

A broader notion of isomorphism can also be stated, where $\alpha$ is allowed to be semi-linear. In view of that we need a few preliminaries on automorphisms of $K$. 

We say that an automorphism $\rho$ of $K$ {\em stabilizes} a given admissible pair $(\sigma,\varepsilon)$ if $\rho\sigma = \sigma\rho$ and $\varepsilon^\rho = \varepsilon$. 

Let $\rho \in \mathrm{Aut}(K)$ stabilize $(\sigma,\varepsilon)$. Then $\rho$ stabilizes both $K_{\sigma,\varepsilon}$ and $K^{\sigma,\varepsilon}$. Thus $\rho$ induces on the group $\overline{K} = K/K_{\sigma,\varepsilon}$ an automorphism $\bar{\rho}$ stabilizing $\overline{K}^\circ = K^{\sigma,\varepsilon}/K_{\sigma,\varepsilon}$. Moreover, $(\bar{t}\circ\lambda)^{\bar{\rho}} = \bar{t}^{\bar{\rho}}\circ\lambda^\rho$ for every element $\bar{t}\in \overline{K}$ and every scalar $\lambda\in K$. Hence the automorphism of $\overline{K}^\circ$ induced by $\bar{\rho}$ is a bijective $\rho$-semi-linear mapping of the $K$-vector space $\overline{K}^\circ$.  

Given a $\circ$-closed subgroup $\overline{R}$ of $\overline{K}$, let $\overline{R}^{\bar{\rho}}$ be the image of $\overline{R}$ by $\bar{\rho}$. Then $\overline{R}^{\bar{\rho}}$ is $\circ$-closed and $\bar{\rho}$ induces an isomorphism from $\overline{K}/\overline{R}$ to $\overline{K}/\overline{R}^{\bar{\rho}}$. Clearly, for every element $\bar{t}+\overline{R}$ of $\overline{K}/\overline{R}$ and every $\lambda\in K$ we have  
\[((\bar{t}+\overline{R})\circ\lambda)^{\bar{\rho}} = (\bar{t}^{\bar{\rho}} + \overline{R}^{\bar{\rho}})\circ\lambda^\rho = (\bar{t}+\overline{R})^{\bar{\rho}}\circ\lambda^\rho.\]
We can now loose our previous definition of isomorphism.   

Let $\overline{R}$ and $\overline{R}'$ be two $\circ$-closed subgroups of $K$. We say that two generalized $(\sigma,\varepsilon)$-quadratic forms $q:V\rightarrow \overline{K}/\overline{R}$ and $q':V'\rightarrow\overline{K}/\overline{R}'$ are {\em weakly isomorphic} if there exists an automorphism $\rho$ of $K$ stabilizing $(\sigma,\varepsilon)$ and such that $\overline{R}^{\bar{\rho}} = \overline{R}'$ and a $\rho$-semi-linear mapping $\alpha:V\rightarrow V'$ such that 
$q'(\alpha(x)) = q(x)^{\bar{\rho}}$ for every $x\in V$. 

\subsection{Proportionality}\label{proportional gen quadric sec}

For $i = 1, 2$ let $(\sigma_i,\varepsilon_i)$ be an admissible pair of $K$ and $\overline{R}_i$ a $\circ_{\sigma_i}$-closed subgroup of $\overline{K}^{\sigma_i,\varepsilon_i} = K/K_{\sigma_i,\varepsilon_i}$ (notation as in Subsection \ref{proportional pairs sec}). Let $q_i:V\rightarrow \overline{K}^{\sigma_i,\varepsilon_i}/\overline{R}_i$ be a non-trivial generalized $(\sigma_i,\varepsilon_i)$-quadratic form and let $f_i$ be its sesquilinearization. We say that $q_1$ and $q_2$ are {\em proportional} if there exists a scalar $\kappa\in K-\{0\}$ such that $(\sigma_2,\varepsilon_2) = \kappa\cdot(\sigma_1,\varepsilon_1)$, $\overline{R}_2 = \kappa \overline{R}_1$ and $q_2(x) = \kappa q_1(x)$ for every $x\in V$. If this is the case then we write $q_2 = \kappa q_1$. Clearly, if $q_2 = \kappa q_1$ then $f_2 = \kappa f_1$ and $S_{q_1}  = S_{q_2}$.

\begin{theo}\label{proportional gen quadric theo}
Let $q_1:V\rightarrow \overline{K}^{\sigma_1,\varepsilon_1}/\overline{R}_1$ and $q_2:V\rightarrow\overline{K}^{\sigma_2,\varepsilon_2}/\overline{R}_2$ be generalized pseudo-quadratic forms such that $S_{q_1} = S_{q_2}$. Assume that the polar space $S := S_{q_1} = S_{q_2}$ has non-degenerate rank at least $2$. Then $q_1$ and $q_2$ are proportional.
\end{theo}
\pr By the same argument used in the proof of Proposition \ref{proportional quadric prop} we obtain that $f_1$ and $f_2$ are proportional. Thus, modulo replacing $q_1$ with $\kappa q_1$ for a suitable $\kappa \in K-\{0\}$ me may assume that $f_1 = f_2 = f$, say. Hence $(\sigma_1,\varepsilon_1) = (\sigma_2,\varepsilon_2)$ and $\overline{K}^{\sigma_1,\varepsilon_1} = \overline{K}^{\sigma_2,\varepsilon_2} =: \overline{K}$.  We must prove that we also have $q_1 = q_2$. As $f_1 = f_2 = f$, we can choose the same facilitating form $g$ for $q_1$ and $q_2$, defining it as in
(\ref{facilitating-E1}) of Subsection \ref{facilitating form}. So, for every $x\in V$, we can choose the same representative $\bar{t}_x \in \overline{K}$ for both $q_1(x)$ and $q_2(x)$. In order to prove that $q_1 = q_2$ we must only show that $\overline{R}_1 = \overline{R}_2$. 

Let $\bar{r}\in \overline{R}_1$. Let $a$ and $b$ be two vectors such that $f(a,b) = 1$ and $[a], [b] \in S$ ($:= S_{q_1} = S_{q_2}$). Such a pair of vectors exists in view of the hypotheses made on $S$. Let $r\in K$ be such that $\bar{r}\in \overline{R}_1$. Then $q_1(a+br) = \bar{r}+\overline{R}_1 = \overline{R}_1$. Hence $[a+br]\in S$. On the other hand, $q_2(a+br) = \bar{r}+\overline{R}_2$. As $[a+br]\in S$, the vector $a+br$ is also $q_2$-singular, namely $\bar{r}\in \overline{R}_2$. It follows that $\overline{R}_1\subseteq\overline{R}_2$. By symmetry, $\overline{R}_2\subseteq \overline{R}_1$. Hence $\overline{R}_1 = \overline{R}_2$.  \eop 

\section{Quotients and covers}

In this section $q:V\rightarrow \overline{K}/\overline{R}$ is a given non-trivial generalized $(\sigma,\varepsilon)$-quadratic form, $f:V\times V\rightarrow K$ is its sesquilinearization and $S_q = (P_q, L_q)$ is the polar space associated to $q$. As $q$ is non-trivial, the form $f$ is non-trivial as well, by Proposition \ref{triviality}. Moreover, $\overline{R}$ is a vector subspace of $\overline{K}^\circ$, by Theorem \ref{forms-theo}, $(1)$.  

We assume that $P_q$ is not totally singular. Hence it spans $\mathrm{PG}(V)$ (Proposition \ref{span}). Therefore the inclusion mapping $e_q:S_q\rightarrow\mathrm{PG}(V)$ is an embedding of $S_q$ in $\mathrm{PG}(V)$. 

Recall that $[\mathrm{Rad}(q)] = [\mathrm{Rad}(f)]\cap P_q$ is the radical of $S_q$. 

\subsection{Quotients}\label{quotient sec}

According to the definitions stated in Subsection \ref{Emb}, a subspace $U$ of $V$ defines a quotient of the embedding $e_q:S_q\rightarrow \mathrm{PG}(V)$ precisely when $[U]\cap P_q = \emptyset$ and $[U]\cap[a,b] = \emptyset$ for any two distinct points $[a], [b]\in P_q$. 

\begin{prop}\label{quot1}
A subspace $U$ of $V$ defines a quotient of the embedding $e_q$ if and only if $U\subseteq \mathrm{Rad}(f)$ and $U\cap \mathrm{Rad}(q) = 0$. 
\end{prop}
\pr This proposition is a special case of the following more general statement on quotients of embeddings of point-line geometries.

Let $e:G\rightarrow \mathrm{PG}(V)$ be a projective embedding of a point-line geometry $G = (P, L)$. Let $W$ be a subspace of $V$ such that a point $[v]$ of $\mathrm{PG}(V)-e(P)$ belongs to $[W]$ if and only if every line of $\mathrm{PG}(V)$ through $[v]$ meets $e(P)$ in at most one point. Then a subspace $U$ of $V$ defines a quotient of the embedding $e$ if and only if $[U]\cap e(P) = \emptyset$ and $U \subseteq W$.

The proof of this claim is easy. We leave it to the reader. In view of the above, in order to prove Proposition \ref{quot1} we only must prove that a point $[v]$ of $\mathrm{PG}(V)- P_q$ belongs to $[\mathrm{Rad}(f)]$ if and only if every projective line through $[v]$ meets $P_q$ in at most one point.

Given a point $[v]\not \in P_q$, assume firstly that every projective line through $[v]$ meets $P_q$ in at most one point. Let $[a]\in P_q$. Then $q(a) = 0_{\bar{R}}$. Consequently, $q(a\lambda+v) = q(v)+ (\overline{\lambda^\sigma f(a,v)} + \overline{R})$ for any $\lambda\in K$. It follows that if $f(a,v) \neq 0$ then a scalar $\lambda\in K$ exists such that  $q(a\lambda + v) = 0_{\bar{R}}$. If this is the case then $[a,v]$ meets $P_q$ in at least two points, namely $[a]$ and $[a\lambda+v]$, a contradiction with the hypotheses made on $[v]$. Therefore $f(a,v) = 0$. As this holds for any $[a]\in P_q$, we obtain that $P_q\subseteq [v^\perp]$. However $P_q$ spans $\mathrm{PG}(V)$, by assumption. Hence $V = v^\perp$, namely $v\in \mathrm{Rad}(f)$.

Conversely, let $v\in\mathrm{Rad}(f)$. Let $[a]\in P_q$. Then $q(a) = 0_{\bar{R}}$ and $f(a,v) = 0$ while $q(v)\neq 0_{\bar{R}}$ as $[v]\not\in P_q$ by assumption. Hence $q(a\lambda + v) = q(v) \neq 0_{\bar{R}}$ for any $\lambda\in K$. This shows that $[a, v]\cap P_q = \{[a]\}$. Therefore every projective line through $[v]$ meets $P_q$ in at most one point.   \eop  

\bigskip

The next corollary immediately follows from Proposition \ref{quot1}.

\begin{co}\label{quot2}
If $\mathrm{Rad}(q) = \mathrm{Rad}(f)$ then the embedding $e_q$ does not admit any proper quotient.
\end{co}

For the rest of this subsection we assume that $\mathrm{Rad}(q)\neq \mathrm{Rad}(f)$. Hence $S_q$ is a proper subspace of $S_f$. Consequently, $(\sigma,\varepsilon)$ is not of trace type. In particular, $\mathrm{char}(K) = 2$. 

Let $U$ be a subspace of $\mathrm{Rad}(f)$ with $U\cap \mathrm{Rad}(q) = 0$. By Proposition \ref{hom}, the restriction of $q$ to $U$ is an injective linear mapping from $U$ to the $K$-vector space $\overline{K}^\circ/\overline{R}$. Hence the image $q(U)$ of $U$ by $q$ is a vector subspace of $\overline{K}^\circ/\overline{R}$. Therefore there exists a unique subspace $\overline{R}_U$ of $\overline{K}^\circ$ containing $\overline{R}$ and such that $\overline{R}_U/\overline{R} = q(U)$. Let $q_U:V/U\rightarrow \overline{K}/\overline{R}_U$ be the mapping defined as follows:
\[q_U(x+U) = \bar{t}+\overline{R}_U ~ \mbox{for an element}~ t\in K ~ \mbox{such that} ~ \bar{t}+\overline{R} = q(x).\]
\begin{lemma}\label{quot4}
The mapping $q_U$ is well defined.
\end{lemma}
\pr Clearly, the coset $\bar{t}+\overline{R}_U$ does not depend on the choice of the representative $\bar{t}$ of $q(x)$. It remains to prove that it neither depends on the choice of the vector $x$ in the coset $x+U$. 

Given $u\in U$, let $x' = x+u$ and let $\bar{t}'$ be a representative of $q(x')$. Then $q(x') = q(x+u) = q(x)+q(u)+f(x,u) = q(x) + q(u)$ because $u\in U \subseteq \mathrm{Rad}(f)$. However $q(u)\in \overline{R}_U/\overline{R}$ by definition of $\overline{R}_U$. Therefore $\bar{t}-\bar{t}'\in \overline{R}_U$, namely $\bar{t}+\overline{R}_U = \bar{t}'+\overline{R}_U$.  \eop 

\bigskip

The sesquilinearization $f$ of $q$ induces a trace-valued $(\sigma,\varepsilon)$-sesquilinear form $f_U$ on $V/U$. Explicitly,
\[f_U(x+U,y+U) := f(x,y).\]
This definition is consistent. Indeed, since $U\subseteq \mathrm{Rad}(f)$, we have $f(x+u,y+v) = f(x,y)$ for any choice of $u, v\in U$. It is clear that, since $f$ is trace-valued and non-trivial, $f_U$ is trace-valued and non-trivial as well. 

The proof of the following lemma is straightforward. We leave it to the reader.  

\begin{lemma}\label{quot5}
The mapping $q_U$ is a generalized $(\sigma,\varepsilon)$-quadratic form. The form $f_U$ induced by $f$ on $V/U$ is a sesquilinearization of $q_U$.  
\end{lemma}

As $f_U$ is non-trivial, the form $q_U$ is non-trivial if and only if $\overline{R}_U\neq \overline{K}$, by Proposition \ref{triviality}. If this is the case then $f_U$ is the unique sesquilinearization of $q_U$, by
Lemma \ref{form}. Finally, Lemma \ref{quot5} and claim (1) of Theorem \ref{forms-theo} imply the following:

\begin{co}\label{quot5-bis}
Let $q_U$ be non-trivial. Then $\overline{R}_U\subseteq \overline{K}^\circ$.
\end{co}

We call $q_U$ the {\em quotient} of $q$ by $U$. According to the notation of Subsection \ref{Sq gen quadr}, when $q_U$ is non-trivial we denote by $P_{q_U}$ and $L_{q_U}$ the set of $q_U$-singular points and totally $q_U$-singular lines of $\mathrm{PG}(V/U)$, respectively. So, $S_{q_U} = (P_{q_U},L_{q_U})$ is the polar space associated to $q_U$ in $\mathrm{PG}(V/U)$.   

\begin{theo}\label{quot6}
Let $\pi_U:V\rightarrow V/U$ be the projection of $V$ onto $V/U$. 

\begin{itemize}
\item[$(1)$] Let $q_U$ be non-trivial. Then $\pi_U$ induces an isomorphism from $S_q$ to $S_{q_U}$.
\item[$(2)$] Let $q_U$ be trivial. Then both forms $f$ and $f_U$ are alternating and $\pi_U$ induces an isomorphism from $S_q$ to the polar space $S_{f_U}$ associated to $f_U$.  
\end{itemize}
\end{theo}
\pr As $U$ defines a quotient of $S_q$, every coset $x+U$ of $U$ in $V$ contains at most one $q$-singular vector. Therefore $\pi_U$ induces and injective mapping on $P_q$. We firstly prove the following:

\begin{itemize}
\item[$(*)$] For every non-zero vector $x\in V$ we have $q_U(x+U) = 0_{\bar{R}_U}$ if and only if $x+u$ is $q$-singular for some $u\in U$.
\end{itemize}
The coset $x+U$ contains a $q$-singular vector if and only if $q(x+u)\in \overline{R}$ for some vector $u\in U$, namely $q(x)+q(u)\in \overline{R}$. (Recall that $f(x,u) = 0$ since $U\subseteq \mathrm{Rad}(f)$). If this is the case then $q(x) \in \overline{R}_U/\overline{R}$, namely $q_U(x+U) = 0_{\bar{R}_U}$. Conversely,  let $q_U(x+U) = 0_{\bar{R}_U}$. Then there exists an element $\bar{t}\in\overline{R}_U$ such that $q(x) = \bar{t}+\overline{R}$. By definition of $\overline{R}_U$, we have $\bar{t}+\overline{R} = q(u)$ for some $u\in U$. Hence $q(x-u) = 0_{\bar{R}}$, namely $x-u$ is $q$-singular.  Claim $(*)$ is proved. 

Let $q_U$ be non-trivial. By $(*)$, the projection $\pi_U$ induces a bijection from $P_q$ to $P_{q_U}$. Two $q_U$-singular points $[x+U]$ and $[y+U]$ of $\mathrm{PG}(V/U)$ are collinear in $S_{q_U}$ if and only if $f_U(x+U,y+U) = 0$. By the definition of $f_U$, this condition is equivalent to $f(x,y) = 0$, which in its turn characterizes the collinearity of $[x]$ and $[y]$. Claim (1) of the theorem is proved.

Let $q_U$ be trivial. Then $(*)$ shows that $\pi_U$ induces a bijection from $P_q$ to the set of points of $\mathrm{PG}(V/U)$. In other words, every coset $x+U$ of $U$ other than $U$ contains exactly one $q$-singular vector. We may assume that in a symbol as $x+U$ the letter $x$ stands for the unique $q$-singular vector of $x+U$. With this convention, $f_U(x+U,x+U) = f(x, x)$ (by definition of $f_U$) and $f(x,x) = 0$ because $x$ is $q$-singular, whence $f$-isotropic. It follows that $f_U(x+U,x+U) = 0$ for every coset $x+U$. Thus, $f_U$ is alternating. Moreover, for any vector $x\in V$ we have $f(x,x) = f_U(x+U,x+U)$ by definition of $f_U$ and $f_U(x+U,x+U) = 0$ since $f_U$ is alternating. Hence $f(x,x) = 0$ for every $x\in V$, namely $f$ is alternating as well. Turning to $S_q$, two points $[x], [y]\in S_q$ are collinear in $S_q$ if and only if $f(x,y) = 0$, equivalently $f_U(x+U,y+U) = 0$, namely $x+U$ and $y+U$ represent collinear points of $S_{f_U}$. Therefore $\pi_U$ maps $S_q$ isomorphically onto $S_{f_U}$, as claimed in (2).  \eop  

\subsection{Covers}\label{cover-sec} 

\subsubsection{Construction and properties of covers}\label{cover-subsec}
 
Let $\overline{S} \oplus \overline{T} = \overline{R}$ be a direct sum decomposition of the $K$-vector space $\overline{R}$. Put $V^{\overline{S}} := V\oplus\overline{S}$ (direct sum of $K$-vector spaces). Define $f^{\overline{S}}:V^{\overline{S}}\times V^{\overline{S}}\rightarrow K$ as follows:
\[f^{\overline{S}}(x+\bar{r}, y+\bar{s}) = f(x,y) ~~ \mbox{for all}~ x,y\in V ~\mbox{and}~ \bar{r}, \bar{s}\in\overline{S}.\]
It is easy to see that $f^{\overline{S}}$ is a trace-valued $(\sigma,\varepsilon)$-sesquilinear form with $\mathrm{Rad}(f^{\overline{S}}) = \mathrm{Rad}(f)\oplus\overline{S}$. Clearly, $f$ is isomorphic to the form induced by $f^{\overline{S}}$ on $V^{\overline{S}}/\overline{S}$ ($\cong V$). 

Let $E = (e_i)_{i\in I}$ be a $q$-singular basis of $V$ and let $g_E$ be the facilitating form associated to $E$ (see definition (\ref{facilitating-E1}) in Subsection \ref{facilitating form}). We define a mapping $q_E^{\overline{S},\overline{T}}:V^{\overline{S}}\rightarrow \overline{K}/\overline{T}$ as follows:
\[q_E^{\overline{S},\overline{T}}(x+\bar{r}) = \overline{g_E(x,x)} + \bar{r} + \overline{T} ~~\mbox{for any}~ x\in V ~\mbox{and any}~ \bar{r}\in \overline{S}. \]
In particular, $q_E^{\overline{S},\overline{T}}(x) = \overline{g_E(x,x)} + \overline{T}$ and $q_E^{\overline{S},\overline{T}}(\bar{r}) = \bar{r}+\overline{T}$. 

\begin{theo}\label{extension-theo}
The mapping $q_E^{\overline{S},\overline{T}}$ is a non-trivial generalized $(\sigma,\varepsilon)$-quadratic form and $f^{\overline{S}}$ is its sesquilinearization. 
\end{theo}
\pr Let $x = \sum_i e_i\lambda_i$ and $\bar{r}\in\overline{S}$. According to the definition of $q_E^{\overline{S},\overline{T}}$ we have 
\[q_E^{\overline{S},\overline{T}}((x+\bar{r})\lambda) = q_E^{\overline{S},\overline{T}}(x\lambda + \bar{r}\circ\lambda) = \sum_{i<j}\overline{\lambda^\sigma\lambda_i^\sigma f(e_i,e_j)\lambda_j\lambda} + \bar{r}\circ\lambda + \overline{T} = \]
\[=  (\sum_{i<j}\overline{\lambda_i^\sigma f(e_i,e_j)\lambda_j})\circ\lambda + \bar{r}\circ\lambda + \overline{T} =  q_E^{\overline{S},\overline{T}}(x+\bar{r})\circ\lambda.\]
So, $q_E^{\overline{S},\overline{T}}$ satisfies condition $(Q'1)$. Turning to $(Q'2)$, let $x = \sum_ie_i\lambda_i$, $y = \sum_ie_i\mu_i$ and $\bar{r}, \bar{s}\in\overline{S}$. Then 
\begin{equation}\label{extension-eq1}
\begin{array}{l}
q_E^{\overline{S},\overline{T}}((x+\bar{r})+(y+\bar{s})) = q_E^{\overline{S},\overline{T}}((x+y)+(\bar{r}+\bar{s})) = \\
{}\\
= \sum_{i<j}\overline{f(e_i,e_j)}\circ(\lambda_j+\mu_j) + \bar{r} + \bar{s} + \overline{T}.
\end{array}
\end{equation}
On the other hand,
\begin{equation}\label{extension-eq2}
\begin{array}{l}
q_E^{\overline{S},\overline{T}}(x+\bar{r}) + q_E^{\overline{S},\overline{T}}(y+\bar{s}) = \\
{}\\
=  \sum_{i<j}\overline{ f(e_i,e_j)}\lambda_j + \sum_{i<j}\overline{ f(e_i,e_j)}\mu_j + \bar{r} + \bar{s} + \overline{T}.
\end{array}
\end{equation}
Moreover, 
\begin{equation}\label{extension-eq3}
f^{\overline{S}}(x+\bar{r},y+\bar{s}) = f(x,y) = \sum_{i<j}(\lambda_i^\sigma f(e_i,e_j)\mu_j.
\end{equation}
By (\ref{extension-eq1}), (\ref{extension-eq2}) and (\ref{extension-eq3}) and recalling that 
\[\begin{array}{l}
\mu_i^\sigma f(e_i,e_j)\lambda_j - \lambda_j^\sigma f(e_j,e_i)\mu_i  =  \mu_i^\sigma f(e_i,e_j)\lambda_j - \lambda_j^\sigma f(e_i,e_j)^\sigma\varepsilon\mu_i  = \\
 = \mu_i^\sigma f(e_i,e_j)\lambda_j - (\mu_i^\sigma f(e_i,e_j)\lambda_i)^\sigma\varepsilon ~ \in K_{\sigma,\varepsilon}
\end{array}\]
we obtain 
$$q_E^{\overline{S},\overline{T}}((x+\bar{r})+(y+\bar{s}))-q_E^{\overline{S},\overline{T}}(x+\bar{r})-q_E^{\overline{S},\overline{T}}(y+\bar{s}) - (\overline{f(x+\bar{r},y+\bar{s})}+\overline{T}) = $$
\[= \sum_{i<j}(\overline{\lambda_j^\sigma f(e_j,e_i)\mu_i} +\overline{\lambda_i^\sigma f(e_i,e_j)\mu_j} -\sum_{i,j}\overline{\lambda^\sigma f(e_i,e_j)\mu_j} + \overline{T} = \]
\[= \sum_i\overline{\lambda_i^\sigma f(e_i,e_i)\mu_i} + \overline{T} = \overline{T}.\]
(Recall that $f(e_i,e_i) = 0$ since $q(e_i) = 0_{\overline{R}}$ by assumption.) Finally,
 $$q_E^{\overline{S},\overline{T}}((x+\bar{r})+(y+\bar{s}))-q_E^{\overline{S},\overline{T}}(x+\bar{r})-q_E^{\overline{S},\overline{T}}(y+\bar{s}) - (\overline{f(x+\bar{r},y+\bar{s})} + \overline{T}) = \overline{T}.$$
Property $(Q'2)$ is proved. The non-triviality of $q_E^{\overline{S},\overline{T}}$ immediately follows from the fact that $q$ is non-trivial by assumption.   \eop 

\bigskip

We say that $q_E^{\overline{S},\overline{T}}$ is the {\em cover} of $q$ via $(\overline{S},\overline{T}$) {\em based} at $E$ (a {\em cover} of $q$, for short). A motivation for this definition is given by the following theorem.  

\begin{theo}\label{covers-main}
The subspace $\overline{S}$ of $V^{\overline{S}}$ defines a quotient $(q_E^{\overline{S},\overline{T}})_{\overline S}$ of $q_E^{\overline{S},\overline{T}}$. With an obvious identification of $V^{\overline{S}}/\overline{S}$ with $V$, we have $(q_E^{\overline{S},\overline{T}})_{\overline{S}} = q$. 
\end{theo}

The proof is straightforward. We leave it to the reader.  By combining this theorem with Theorem \ref{quot6} we immediately obtain the following: 

\begin{co}\label{cover-main-cor}
The polar space associated to $q^{\overline{S},\overline{T}}_E$ in $\mathrm{PG}(V^{\overline{S}})$ is isomorphic to the polar space $S_q$ associated to $q$.
\end{co}

Theorem \ref{covers-main} can be rephrased in the language of embeddings, but in view of that we need a few more definitions. For $\bar{r}\in \overline{R}$, let $\theta(\bar{r})$ be the projection of $\bar{r}$ onto $\overline{S}$ along $\overline{T}$, namely $\theta(\bar{r})$ is the unique element of $\overline{S}\cap(\bar{r}+\overline{T})$. For every $q$-singular vector $x\in V$, the subspace $\langle x, \overline{S}\rangle$ of $V^{\overline{S}}$ contains a unique $q_E^{\overline{S},\overline{T}}$-singular point, represented by the vector $x-\theta(\overline{g_E(x,x)})$. Put
\begin{equation}\label{lifting of eq}
e_{q, E}^{\overline{S},\overline{T}}([x]) := [x-\theta(\overline{g_E(x,x)})].
\end{equation}
The following is straightforward. We leave its proof to the reader. 
\begin{theo}\label{covers-main-bis}
The mapping $e_{q,E}^{\overline{S},\overline{T}}$ is a projective embedding of $S_q$ in $\mathrm{PG}(V^{\overline{S}})$. The image $e_{q, E}^{\overline{S},\overline{T}}(S_q)$ of $S_q$ by $e^{\overline{S},\overline{T}}_{q,E}$ is the polar space associated to $q^{\overline{S},\overline{T}}_E$ in $\mathrm{PG}(V^{\overline{S}})$. Moreover, if $\pi_{\overline{S}}$ is the projection of $V^{\overline{S}}$ onto $V^{\overline{S}}/\overline{S}$, then the canonical isomorphism from $V^{\overline{S}}/\overline{S}$ to $V$ yields an isomorphism from the composition $\pi_{\overline{S}}\cdot e_{q,E}^{\overline{S},\overline{T}}$ to the inclusion embedding $e_q:S_q\rightarrow \mathrm{PG}(V)$.
\end{theo}

We call $e_{q,E}^{\overline{S},\overline{T}}$ the {\em lifting} of $e_q$ to $V^{\overline{S}}$ {\em based} at $E$. 

\bigskip

\noindent
{\bf Remark.} We have assumed that $q$ is non-trivial since the very beginning of Section 4, however the previous construction can be repeated when $q$ is trivial. In that case we choose a sesquilinearization $f$ of $q$ and we define $q_E^{\overline{S},\overline{T}}$ with the help of $f$, as in the non-trivial case, but the form $q^{\overline{S},\overline{T}}_E$ now depends on the particular choice of $f$. The form $q_E^{\overline{S},\overline{T}}$ is non-trivial provided that $\overline{S}\neq \{\bar{0}\}$. It is still true that $q$ is a quotient of $q_E^{\overline{S},\overline{T}}$, but Corollary \ref{cover-main-cor} must be rephrased as follows: the polar space associated to $q^{\overline{S},\overline{T}}_E$ in $\mathrm{PG}(V^{\overline{S}})$ is isomorphic to $S_f$ (compare Theorem \ref{quot6}, (2)).    

\subsubsection{Independence of $q_E^{\overline{S},\overline{T}}$ from the choice of $E$}

Our definition of $q_E^{\overline{S},\overline{T}}$ rests on the choice of a particular ordered $q$-singular basis $E$. In this subsection we shall prove that this choice is ultimately irrelevant: different choices lead to isomorphic forms. 

Given two $q$-singular bases $E$ and $E'$, let $\delta_{E,E'}$ be the difference-map of the pair $(E,E')$ (see Subsection \ref{facilitating form}). Recall that $\delta_{E,E'}(x)\in\overline{R}_{E,E'}\subseteq \overline{R}$, by Lemma \ref{irrelevant}. Hence $\theta(\delta_{E,E'}(x))$ is defined for every $x\in V$, where $\theta$ is the projection of $\overline{R}$ onto $\overline{S}$ along $\overline{T}$, as in (\ref{lifting of eq}). In view of the definition of $\delta_{E,E'}$, the following holds for every vector $x\in V$:
\[x-\theta(\overline{g_{E'}(x,x)}) = x-\theta(\overline{g_E(x,x)}) + \theta(\delta_{E,E'}(x)).\]
Let $\Delta_{E,E'}:V^{\overline{S}}\rightarrow V^{\overline{S}}$ be the mapping defined as follows:
\[\Delta_{E,E'}(x+\bar{r}) = x+\theta(\delta_{E,E'}(x)) + \bar{r} ~ \mbox{for any} ~ x\in V ~\mbox{and} ~\bar{r}\in \overline{S}\]

\begin{theo}\label{irrelevant-meglio}
The mapping $\Delta_{E,E'}$ is linear and bijective, it fixes $\overline{S}$ elementwise and yields an isomorphism from $q_E^{\overline{S},\overline{T}}$ to $q_{E'}^{\overline{S},\overline{T}}$. Explicitly, 
\begin{equation}\label{irrelevant-eq}
q_E^{\overline{S},\overline{T}}(x+\bar{r}) = q_{E'}^{\overline{S},\overline{T}}(\Delta_{E,E'}(x+\bar{r}))
\end{equation}
for any $x\in V$ and $\bar{r}\in \overline{S}$. Consequently, $\Delta_{E,E'}$ is an isomorphism of embeddings from the lifting $e_{q,E}^{\overline{S},\overline{T}}$ of $e_q$ based at $E$ to the lifting $e_{q,E'}^{\overline{S},\overline{T}}$ of $e_q$ based at $E'$.  
\end{theo}
\pr By Lemma \ref{irrelevant}, the difference-map $\delta_{E,E'}$ is a linear mapping from $V$ to $\overline{R}_{E,E}$. Hence $\Delta_{E,E'}$ is linear. Clearly, $\Delta_{E,E'}$ fixes $\overline{S}$ elementwise. Moreover the composition of $\Delta_{E,E'}$ with the projection of $V^{\overline{S}}$ onto $V$ along $\overline{S}$ induces the identity mapping on $V$. Therefore $\Delta_{E,E'}$ is bijective.  We have 
\[q_E^{\overline{S},\overline{T}}(x+\bar{r}) = \overline{g_E(x,x)}+\bar{r}+\overline{T} = \]
\[= \overline{g_{E'}(x,x)} + (\overline{g_E(x,x)}-\overline{g_{E'}(x,x)}) +\bar{r}+\overline{T} =\]
\[= \overline{g_{E'}(x,x)} + \delta_{E,E'}(x) +\bar{r}+\overline{T} = \]
\[= \overline{g_{E'}(x,x)} + \theta(\delta_{E,E'}(x)) +\bar{r}+\overline{T} = q_{E'}^{\overline{S},\overline{T}}.\]
(Recall that $\delta_{E,E'}(x) +\overline{T} = \theta(\delta_{E,E'}(x)) + \overline{T}$, by the definition of $\theta$.) Equation (\ref{irrelevant-eq}) is proved. Exploiting  (\ref{irrelevant-eq}), it is not difficult to prove that $\Delta_{E,E'}$ is an isomorphism from $e_{q,E}^{\overline{S},\overline{T}}$ to $e_{q,E'}^{\overline{S},\overline{T}}$.  \eop

\bigskip

Theorem \ref{irrelevant-meglio} allows us to drop the index $E$ in our notations, thus writing $q^{\overline{S},\overline{T}}$ and $e_q^{\overline{S},\overline{T}}$ for $q_E^{\overline{S},\overline{T}}$ and $e_{q,E}^{\overline{S},\overline{T}}$ whenever the particular choice of the basis $E$ is irrelevant for what we are saying. Accordingly, we may call $q^{\overline{S},\overline{T}}$ and $e_q^{\overline{S},\overline{T}}$ the {\em cover} of $q$ via $(\overline{S},\overline{T})$ and the {\em lifting} of $e_q$ to $V^{\overline{S}}$ respectively, with no mention of the basis $E$. 

\subsubsection{Dominant covers}\label{dominant covers}

As $\overline{S}\oplus\overline{T} = \overline{R}$, we have $\overline{S} = \overline{R}$ if and only if $\overline{T} = \{\bar{0}\}$. When $\overline{T} = \{\bar{0}\}$ the form $q^{\overline{S},\overline{T}} = q^{{\overline{R}},\{\bar{0}\}}$ is pseudo-quadratic with defect equal to $\mathrm{Rad}(f)\oplus \overline{R}$.  

Improper covers are allowed too. We get them by taking $\overline{S} = \{\bar{0}\}$ (whence $\overline{T} = \overline{R}$). Clearly, $q^{\{\bar{0}\},\overline{R}} = q$. 

Notice that we have not assumed that $\overline{R} \neq \{\bar{0}\}$. Indeed the construction of $q^{\overline{S},\overline{T}}$ makes sense even if $\overline{R} = \{\bar{0}\}$, namely $q$ is pseudo-quadratic. In this case $\overline{S} = \overline{T} = \{\bar{0}\}$, hence $q^{\overline{S},\overline{T}} = q$, namely $q$ does not admit any proper cover. Conversely, if $q$ does not admit any proper cover then $\overline{R} = \{\bar{0}\}$. 

We say that $q$ is {\em dominant} if it does not admit any proper cover. By the above, $q$ is dominant if and only if it is pseudo-quadratic. So, the form $q^{\overline{S},\overline{T}}$ is dominant if and only if $\overline{T} = \{\bar{0}\}$. We call $q^{\overline{R},\{\overline{0}\}}$ the {\em dominant cover} of $q$. 

\subsubsection{Quotients versus covers} 

According to Theorem \ref{covers-main}, if $\tilde{q}:\widetilde{V}\rightarrow \overline{K}/\overline{T}$ is a cover of $q:V\rightarrow\overline{K}/\overline{R}$ then $q$ is a quotient of $\tilde{q}$. A converse of this statement also holds.

\begin{theo}\label{covers-main-ter}
Given a subspace $\overline{T}$ of $\overline{K}^\circ$ and a generalized $(\sigma,\varepsilon)$-quadratic form $\tilde{q}:\widetilde{V}\rightarrow \overline{K}/\overline{T}$, let $U$ be a subspace of $\widetilde{V}$ defining a quotient of $\tilde{e}$. Then $\tilde{q}$ is isomorphic to a cover of the quotient $\tilde{q}_U$ of $\tilde{q}$ by $U$.
\end{theo}
\pr Put $V := \widetilde{V}/U$ and $q := \tilde{q}_U:V\rightarrow \overline{K}/\overline{R}$, where $\overline{R} := \overline{T}_U$ is the subspace of $\overline{K}^\circ$ such that $\overline{R}/\overline{T} = \tilde{q}(U)$ (see Subsection \ref{quotient sec}). Let $\overline{S}$ be a complement of $\overline{T}$ in the $K$-vector space $\overline{R}$, let $W$ be a complement of $U$ in $\widetilde{V}$, let $\tilde{\pi}_U$ be the projection of $\widetilde{V}$ onto $V = \widetilde{V}/U$ and $\theta$ the projection of $\overline{R}$ onto $\overline{S}$ along $\overline{T}$. Let $\alpha:\widetilde{V}\rightarrow V^{\overline{S}} = V \oplus \overline{S}$ be the linear mapping defined by the following clauses: $\alpha(w) = \tilde{\pi}_U(w)$ for every $w\in W$ and $\alpha(u) = \theta(\tilde{q}(u))$ for every $u\in U$. As the reader can check, $\alpha$ is an isomorphism from $\tilde{q}$ to $q^{\overline{S},\overline{T}}$. \eop

\begin{co}\label{covers-main-ter-co}
Let $q:V\rightarrow \overline{K}/\overline{R}$ be a non-trivial generalized $(\sigma,\varepsilon)$-quadratic form. Given a vector subspace $\overline{T}$ of $\overline{R}$, let $\overline{S}$ and $\overline{S}'$ be two complements of $\overline{T}$ in $\overline{R}$. Then $q^{\overline{S},\overline{T}}\cong q^{\overline{S}',\overline{T}}$.
\end{co}  
\pr The conclusion follows from the proof of Theorem \ref{covers-main-ter} with $V^{\overline{S}'}$, $q^{\overline{S}',\overline{T}}$ and $\overline{S}'$ in the roles of $\widetilde{V}$, $\tilde{q}$ and $U$ respectively, recalling that $q$ is the quotient of $q^{\overline{S}',\overline{T}}$ over $\overline{S}'$ by Theorem \ref{covers-main}.   \eop   

\subsubsection{Partial independence of $q^{\overline{S},\overline{T}}$ from the choice of $\overline{S}$ and $\overline{T}$} 

In general, if $\overline{R} = \overline{S}\oplus\overline{T}$ and $\overline{R} = \overline{S}'\oplus \overline{T}'$ are two decompositions of $\overline{R}$ then $q^{\overline{S},\overline{T}} \not\cong q^{\overline{S}',\overline{T}'}$. However, with a suitable choice of $\overline{T}'$ the forms $q^{\overline{S},\overline{T}}$ and $q^{\overline{S}',\overline{T}'}$ are weakly isomorphic in the sense of Subsection \ref{weak iso}. Explicitly:

\begin{prop}
With $\overline{S}, \overline{T}, \overline{S}'$ and $\overline{T}'$ as above, suppose that $K$ admits an automorphism $\rho$ stabilizing $(\sigma,\varepsilon)$ and such that the automorphism $\bar{\rho}$ of $\overline{K}$ induced by $\rho$ stabilizes $\overline{R}$ and maps $\overline{T}$ onto $\overline{T}'$. Then the forms $q^{\overline{S},\overline{T}}$ and $q^{\overline{S}',\overline{T}'}$ are weakly isomorphic.
\end{prop}
\pr Given a $q$-singular basis $E$ of $V$ let $\rho_E$ be the $\rho$-semi-linear mapping of $V$ that fixes all vetors of $E$ and, for $x\in V$ and $\bar{r}\in \overline{S}$, set $\rho_E(x+\bar{r}) := \rho_E(x) + \bar{r}^{\bar{\rho}}$. Then $\rho_E$ is a bijective $\rho$-semilinear mapping from $V^{\overline{S}}$ to $V^{\overline{S}^{\bar{\rho}}}$ and
we have
\[(q_E^{\overline{S},\overline{T}}(x+\bar{r}))^{\bar{\rho}} = q_E^{\overline{S}^{\bar{\rho}},\overline{T}^{\bar{\rho}}}(\rho_E(x+\bar{r}))\] 
for every vector $x+\bar{r}$ of $V^{\overline{S}}$. Hence $q^{\overline{S},\overline{T}}$ and $q^{\overline{S}^{\bar{\rho}},\overline{T}^{\bar{\rho}}}$ are weakly isomorphic. However $q^{\overline{S}^{\bar{\rho}},\overline{T}^{\bar{\rho}}} \cong q^{\overline{S}',\overline{T}'}$ by Corollary \ref{covers-main-ter-co} and because $\overline{R}^{\bar{\rho}} = \overline{R}$ and $\overline{T}^{\bar{\rho}} = \overline{T}'$ by assumption. Therefore $q^{\overline{S},\overline{T}}$ and $q^{\overline{S}',\overline{T}'}$ are weakly isomorphic.  \eop 

\section{Forms for embedded polar spaces} 

Throughout this section $S = (P, L)$ is a non-degenerate polar space of rank at least $2$ and $e:S\rightarrow \mathrm{PG}(V)$ is a projective embedding.  So, the image $e(S) = (e(P), e(L))$ of $S$ by $e$ is a full subgeometry of $\mathrm{PG}(V)$, it spans $\mathrm{PG}(V)$ and $e(S)\cong S$. 

Let $K$ be the underlying divison ring of $V$. By Theorem \ref{Theo 1}, an admissible pair $(\sigma,\varepsilon)$ of $K$ and a $(\sigma,\varepsilon)$-sesquilinear form $f:V\times V\rightarrow K$ exist such that $e(S)$ is a subspace of the polar space $S_f = (P_f,L_f)$ associated to $f$. Explicitly, 

\begin{itemize}
\item[(E1)] $e(P)\subseteq P_f$ and, for any two points $[x], [y] \in e(P)$, the line $[x,y]$ of $\mathrm{PG}(V)$ belongs $e(L)$ if and only if $f(x,y) = 0$.
\end{itemize}
Property (E1) implies both the following:
\begin{itemize} 
\item[(E2)] For any two points $[x]$ and $[y]$ of $\mathrm{PG}(V)$ with $[y]\in e(P)$, we have $f(x,y) = 0$ if and only if either the line $[x,y]$ belongs to $e(L)$ or $[x,y]\cap e(P) = \{[y]\}$. 
\item[(E3)]  $e(P)\cap [\mathrm{Rad}(f)] = \emptyset$. 
\end{itemize}
As for (E3), recall that $S$ is non-degenerate by assumption while $f$ might be degenerate. By (E1), (E2) and (E3) and recalling that $e(P)$ spans $\mathrm{PG}(V)$, we also obtain the following:
\begin{itemize}
\item[(E4)] A point $[x]$ of $\mathrm{PG}(V)$ belongs to $[\mathrm{Rad}(f)]$ if and only if every line of $\mathrm{PG}(V)$ through $[x]$ meets $e(P)$ in at most one point. 
\end{itemize}
The form $f$ is uniquely determined up to proportionality (Proposition \ref{proportional forms prop2}). Moreover $f$ is trace-valued by (2) of Proposition \ref{trace-prop2}, since $P_f\supseteq e(P)$ and $e(P)$ spans $\mathrm{PG}(V)$ 

Let $E = (e_i)_{i\in I}$ be a basis of $V$ such that $[e_i]\in e(P)$ for any $i\in I$. Such a basis exists since $e(P)$ spans $\mathrm{PG}(V)$. We call $E$ an $e(S)$-{\em basis} of $V$. Given a total ordering $<$ on $I$, let $g_E(x,y)$ be defined as in (\ref{facilitating-E1}) of Subsection \ref{facilitating form} and put
\[\gamma_E(x) := \overline{g_E(x,x)} = \sum_{i<j}\overline{\lambda_i^\sigma f(e_i,e_j)\lambda_j} ~~ \mbox{for every vector}~ x = \sum_ie_i\lambda_i \in V.\]
\begin{lemma}\label{qE0}
The mapping $\gamma_E$ is a (possibly trivial) $(\sigma,\varepsilon)$-quadratic form, $g_E$ is a facilitating form for $\gamma_E$ and $f$ is a sesquilinearization of $\gamma_E$. The form $\gamma_E$ is trivial if and only if $\sigma = \mathrm{id}_K$, $\varepsilon = -1$ and $\mathrm{char}(K) \neq 2$.  
\end{lemma} 
\pr The first three claims of the lemma are obvious (compare Subsection \ref{facilitating forms quadratic}). The last one follows from the second part of (\ref{sigma-epsilon3}) of Subsection \ref{section 2.1}. \eop  

\bigskip

Let $\overline{R}$ be the closed subgroup of $\overline{K}$ generated by the set $\{\gamma_E(x)\}_{[x]\in e(P)}$ and define a mapping $q:V\rightarrow \overline{K}/\overline{R}$ as follows:
\begin{equation}\label{definition q}
q(x) := \gamma_E(x) + \overline{R}.
\end{equation}
The next lemma easily follows from Lemma \ref{qE0} and the definition of $\overline{R}$.

\begin{lemma}\label{qE1} 
The mapping $q$ defined in {\rm (\ref{definition q})} is a (possibly trivial) generalized $(\sigma,\varepsilon)$-quadratic form. If $q$ is non-trivial then $f$ is the sequilinearization of $q$. In this case $e(S)$ is a subspace of the polar space $S_q = (P_q, L_q)$ associated to $q$. 
\end{lemma}
\pr The first two claims of the lemma are straightforward. As for the third one, note firstly that $S_q$ is a subspace of $S_f$ since $f$ is the sesquilinearization of $q$ by Lemma \ref{qE1}. Clearly, $e(P)\subseteq P_q$. Therefore $e(S)$ is a subspace of $S_q$, as both $e(S)$ and $S_q$ are subspaces of $S_f$.  \eop 

\begin{co}\label{qE1-bis}
If $(\sigma,\varepsilon)$ is of trace type then either $\overline{R} = \overline{K}$ or $\overline{R} = \{\bar{0}\}$.
\end{co}
\pr This statement easily follows from Lemma \ref{qE1} and Corollary \ref{forms-cor}.  \eop 

\bigskip 

Note that, while $\gamma_E$ depends on the choice of the ordered basis $E$, neither $\overline{R}$ nor $q$ depend on that choice (see the final remark of Subsection \ref{facilitating form}). 

\begin{co}\label{qE1-trivial1}
The form $q$ is trivial if and only if $\overline{R} = \overline{K}$. If $\gamma_E$ is trivial then $\overline{R} = \overline{K}$ (whence $q$ is also trivial)
\end{co}
\pr The form $f$ is non-trivial, since $e(S)$ is a subspace of $S_f$ and it is non-degenerate. This fact and Proposition \ref{triviality} imply the first claim of the corollary. According to the last claim of Lemma \ref{qE0}, the form $\gamma_E$ is trivial if and only if $f$ is alternating and $\mathrm{char}(K)\neq 2$. If this is the case then $\overline{R} = \overline{K}$. \eop 

\begin{theo}\label{main1}
Either $q$ is trivial or $e(S) = S_q$.
\end{theo}
\pr Suppose that $q$ is non-trivial. By Corollary \ref{qE1-trivial1}, $\overline{R}$ is a proper subgroup of $\overline{K}$. Moreover $e(S)$ is a subspace of $S_q$, by the last claim of Lemma \ref{qE1}. 

Let $\tilde{q} := q^{\overline{R},\{\bar{0}\}}:V\oplus\overline{R}\rightarrow \overline{K}$ be the dominant cover of $q$ based at $E$, let $\tilde{f} := f^{\overline{R}}$ be the sesquilinearization of $\tilde{q}$ and put $\widetilde{V} = V\oplus\overline{R}$. Clearly, if $\overline{R} = \{\bar{0}\}$ (as when $(\sigma,\varepsilon)$ is of trace type) then $\tilde{q} = q$, $\tilde{f} = f$ and $\widetilde{V} = V$. 

The embedding $e:S\rightarrow \mathrm{PG}(V)$ lifts to an embedding $\tilde{e}:S\rightarrow \mathrm{PG}(\widetilde{V})$, obtained as the composition of $e$ with the lifting of the inclusion embedding $e_q:S_q\rightarrow \mathrm{PG}(V)$ to $\widetilde{V}$ (see definition (\ref{lifting of eq}) of Subsection \ref{cover-subsec}). Let $\widehat{V}$ be the subspace of $\widetilde{V}$ spanned by $\tilde{e}(P)$. We shall prove that $\widehat{V} = \widetilde{V}$.

Put $\widehat{R} := \overline{R}\cap\widehat{V}$ and let $\hat{q}$ and $\hat{f}$ be the forms induced by $\tilde{q}$ and $\tilde{f}$ on $\widehat{V}$. Clearly, all points of $\tilde{e}(P)$ are $\hat{q}$-singular. As $\widehat{V}+\overline{R} = \widetilde{V}$, we have $\widehat{V}/\widehat{R}\cong \widetilde{V}/\overline{R}\cong V$ and $\widehat{R}$ defines a quotient $\hat{q}_{\widehat{R}}$ of $\hat{q}$. Via an obvious identification of $V$ with $\widehat{V}/\widehat{R}$, we may assume that $\hat{q}_{\widehat{R}}$ is defined over $V$. Accordingly, all points of $e(P)$ are $\hat{q}_{\widehat{R}}$-singular. It follows that $\gamma_E(x)$ belongs to the co-defect $\widehat{R}$ of $\hat{q}_{\widehat{R}}$, for every point $[x]\in e(P)$. However, $\overline{R}$ is generated by $\{ \gamma(x)\}_{[x]\in e(P)}$. Therefore $\widehat{R} = \overline{R}$. Hence $\overline{R}\subseteq \widehat{V}$. It is now clear that $\widehat{V} = \widetilde{V}$, namely $\tilde{e}(P)$ spans $\mathrm{PG}(\widetilde{V})$. 

Since $e(S)$ is a subspace of $S_q$, the image $\tilde{e}(S)$ of $S$ by $\tilde{e}$ is a subspace of the polar space $S_{\tilde{q}} = (P_{\tilde{q}}, L_{\tilde{q}})$ associated to $\tilde{q}$. The latter is a subspace of the polar space $S_{\tilde{f}} = (P_{\tilde{f}},L_{\tilde{f}})$ associated to $\tilde{f}$. Hence $\tilde{e}(S)$ is also a subspace of $S_{\tilde{f}}$, namely (E1) holds with $\tilde{e}(S)$ and $\tilde{f}$ in the roles of $e(S)$ and $f$ respectively. Consequently, properties (E2), (E3) and (E4) also hold for $\tilde{e}(S)$ and $\tilde{f}$. 

We shall now prove that $e(S) = S_q$. Suppose the contrary, namely $e(P)\subset P_q$. Then we also have $\tilde{e}(P)\subset P_{\tilde{q}}$. Let $[a]\in P_{\tilde{q}}-\tilde{e}(P)$. Suppose firstly that $[a]\not\in[\mathrm{Rad}(\tilde{f})]$. By (E4), there exist two distinct points $[b], [c]\in \tilde{e}(P)$ such that the line $[b,c]$ contains $[a]$. We have $\tilde{f}(a,a) = \tilde{f}(b,c) = \tilde{f}(c,c) = 0$ since all of $[a], [b]$ and $[c]$ belong to $P_{\tilde{f}}$. On the other hand, the line $[b,c]$ does not belong to $\tilde{e}(L)$, since it contains $[a]$ which, by assumption, does not belong to $\tilde{e}(P)$. Then $f(b,c) \neq 0$ by (E1). Since $f(b,b) = f(c,c) = 0$ while $f(b,c)\neq 0$, the form $\tilde{f}$ induces a non-degenerate form on the subspace $\langle b, c\rangle$ of $\widetilde{V}$. Thus we can apply Proposition 10.3.10 of Buekenhout and Cohen \cite{BC}. By claim $(i)$ of that proposition, $P_{\tilde{q}}\cap[b,c]$ is the smallest subset of $S_{\tilde{f}}\cap[b,c]$ containing $[b]$ and $[c]$ and perspective with respect to the polarity $\delta_{\tilde{f},[b,c]}$ defined by $\tilde{f}$ on the line $[b,c]$. However, $[b], [c]\in \tilde{e}(P)\cap[b,c]\subseteq P_{\tilde{q}}\cap[b,c]$ and $\tilde{e}(P)\cap[b,c]$ is also perspective with respect to $\delta_{\tilde{f},[b,c]}$ by Proposition 10.3.4 of Buekenhout and Cohen \cite{BC}. Hence $\tilde{e}(P)\cap[b,c] = P_{\tilde{q}}\cap [b,c]$. In particular, $[a]\in \tilde{e}(P)$, contrary to our choice of $[a]$. Therefore $[a]\in [\mathrm{Rad}(\tilde{f})]$, namely $[a]\in[\mathrm{Rad}(\tilde{q})]$, as $[a]\in P_{\tilde{q}}$. It follows that $P_{\tilde{q}}-\tilde{e}(P)\subseteq [\mathrm{Rad}(\tilde{q})]$. 

Still with $[a] \in  P_{\tilde{q}}-\tilde{e}(P) \subseteq  [\mathrm{Rad}(\tilde{q})]$,  let $[b]\in \tilde{e}(P)$. As both $[b]$ and $[a]$ are $\tilde{q}$-singular and $[a]\in[\mathrm{Rad}(\tilde{q})]$, the line $[a, b]$ belongs to $L_{\tilde{q}}$. Hence it is totally $\tilde{f}$-isotropic. By (E1), if $[a,b]$ contains a point of $\tilde{e}(P)$ different from $[b]$ then it also belongs to $\tilde{e}(L)$, but this contradicts the choice of $[a]\in P_{\tilde{q}}-\tilde{e}(P)$. Therefore $[a,b]\cap \tilde{e}(P) = \{[b]\}$, namely $[a,b]-\{[b]\}\subseteq P_{\tilde{q}}-\tilde{e}(P)$. However $P_{\tilde{q}}-\tilde{e}(P)\subseteq [\mathrm{Rad}(\tilde{q})]$ and the latter is a subspace of $\mathrm{PG}(\widetilde{V})$. It follows that $[a,b]\subseteq[\mathrm{Rad}(\tilde{q})]$. This forces $[b]\in [\mathrm{Rad}(\tilde{q})]\cap \tilde{e}(P)\subseteq [\mathrm{Rad}(\tilde{f})]\cap\tilde{e}(P)$, a contradiction with (E3). We have reached a final contradiction. Therefore $e(S) = S_q$.   \eop 

\bigskip

Let $\overline{R}\neq \overline{K}$. Then both $q$ and $\gamma_E$ are non-trivial (Corollary \ref{qE1-trivial1}). Let $S_{\gamma_E} = (P_{\gamma_E}, L_{\gamma_E})$ be the polar space associated t to $\gamma_E$ in $\mathrm{PG}(V)$. Clearly, $S_{\gamma_E}$ is a subspace of $S_f$.  

\begin{co}\label{qE2}
Let $\overline{R} \neq \overline{K}$. Then $S_{\gamma_E}$ is a subspace of $e(S)$. If moreover $(\sigma,\varepsilon)$ is of trace type then $S_{\gamma_E} = e(S) = S_f$. 
\end{co}
\pr Clearly $S_{\gamma_E}$ is a subgeometry of $S_q$. Moreover both $S_{\gamma_E}$ and $S_q$ are subspaces of $S_f$. Hence $S_{\gamma_E}$ is a subspace of $S_q$. However $S_q = e(S)$ by Theorem \ref{main1}. Therefore $S_{\gamma_E}$ is a subspace of $e(S)$.

Let $(\sigma,\varepsilon)$ be of trace type. Then $S_{\gamma_E} = S_f$ by Proposition \ref{Sq trace type}. Hence $S_{\gamma_E} = e(S) = S_f$, since $S_{\gamma_E}$ is a subspace of $e(S) = S_q$ which in its turn is a subspace of $S_f$.   \eop 

\bigskip

\noindent
{\bf Remark.} When $(\sigma,\varepsilon)$ is not of trace type it can happen that $S_{\gamma_E}$ is a proper subspace of $e(S)$. If that is the case then the space $S_{\gamma_E}$ depends on the choice of the $e(S)$-basis $E$.

\begin{theo}\label{qE4}
Let $\overline{R} = \overline{K}$. Then $f$ is an alternating form and $e(S) = S_f$.  
\end{theo}
\pr As $\overline{R} = \overline{K}$, the group $\overline{K}$ is generated by the elements $\gamma_E(x)$ for $[x]\in e(S)$. However $e(S)$ is a subspace of $S_f$. Hence $\overline{K}$ is also generated by the elements $\gamma_E(x)$ for $x$ such that $f(x,x) = 0$. 

Given $x = \sum_ie_i\lambda_i$, let $t := \sum_{i<j}\lambda^\sigma_i f(e_i,e_j)\lambda_j$. Then
\[f(x,x) = \sum_{i,j}\lambda_i^\sigma f(e_i,e_j)\lambda_j = \sum_{i\neq j}\lambda_i^\sigma f(e_i,e_j)\lambda_j + \sum_i\lambda_i^\sigma f(e_i,e_i)\lambda_i.\]
However $f(e_i,e_i) = 0$ for every $i\in I$ because $[e_i]\in e(P)\subseteq P_f$. Therefore
\[\begin{array}{l}
0 = \sum_{i\neq j}\lambda_i^\sigma f(e_i,e_j)\lambda_j  = \sum_{i<j}\lambda_i^\sigma f(e_i,e_j)\lambda_j + \sum_{i > j}\lambda_i^\sigma f(e_i,e_j)\lambda_j = \\
= \sum_{i<j}\lambda_i^\sigma f(e_i,e_j)\lambda_j + \sum_{j > i}\lambda_j^\sigma f(e_j,e_i)\lambda_i =  \\
= \sum_{i<j}\lambda_i^\sigma f(e_i,e_j)\lambda_j + \sum_{i < j}(\lambda_i^\sigma f(e_i,e_j)\lambda_j)^\sigma\varepsilon = \\
= \sum_{i<j}\lambda_i^\sigma f(e_i,e_j)\lambda_j + (\sum_{i < j}\lambda_i^\sigma f(e_i,e_j)\lambda_j)^\sigma\varepsilon  = t + t^\sigma\varepsilon.
\end{array}\]
Hence $f(x,x) = 0$ if and only if $t = -t^\sigma\varepsilon$, namely $t\in K^{\sigma,\varepsilon}$. However $\overline{K}$ is generated by the values $\gamma_E(x)$ with $f(x,x) = 0$. Therefore $K = K^{\sigma,\varepsilon}$. The latter holds precisely when $\varepsilon = -1$ and $\sigma = \mathrm{id}_K$, by the first claim of (\ref{sigma-epsilon3}) of Subsection \ref{section 2.1}.  

So, $\sigma = \mathrm{id}_K$ and $\varepsilon = -1$. In particular, $K$ is a field. If $\mathrm{char}(K) \neq 2$ then $f$ is alternating. Let $\mathrm{char}(K) = 2$. Then $f$ is a symmetric bilinear form. However, $f$ is also trace-valued. It is well known that the alternating forms are the only trace-valued symmetric bilinear forms in characteristic $2$. Hence $f$ is alternating.   

We still must prove that $e(S) = S_f$. This can be proved with the help of Theorem \ref{Theo 3}, but according to the philosophy we have chosen in this paper, we prefer not to use that theorem.

We firstly assume that $\mathrm{char}(K) \neq 2$. By way of contradiction, suppose that $P_f\not\subseteq e(P)$ and let $[a]\in P_f- e(P)$. Assume that $[a]\not\in [\mathrm{Rad}(f)]$. By (E4), there exists at least one line $l$ of $\mathrm{PG}(V)$ containing $[a]$ and intersecting $e(P)$ in at least two points. By Proposition 10.3.4 of Buekenhout and Cohen \cite{BC}, the set $e(P)\cap l$ is perspective with respect to the polarity $\delta_{f,l}$ defined by $f$ on the line $l$. However, according to Bueknehout and Cohen \cite[Proposition 10.3.10(ii)]{BC}, the line $l$ does not contain any proper subset of size at least two and perspective with respect to $\delta_{f,l}$. Therefore $l = e(P)\cap l$. This contradicts the choice of $[a]\not\in e(P)$. We must conclude that $[a]\in [\mathrm{Rad}(f)]$. So, $P_f-e(P) \subseteq [\mathrm{Rad}(f)]$. With $[a]\in P_f-e(P) \subseteq [\mathrm{Rad}(f)]$, let $[b]\in e(P)$. Then $[a,b]\cap e(P) = \{[b]\}$ by (E1). Consequently $[a,b]-\{[b]\}\subseteq [\mathrm{Rad}(f)]$. However $\mathrm{Rad}(f)$ is a subspace of $V$. Hence $[b]\in [\mathrm{Rad}(f)]$, in contradiction with (E3). Therefore $e(S) = S_f$. 

Let now $\mathrm{char}(K) = 2$. Then $K_{\sigma,\varepsilon} = 0$, $K^{\sigma,\varepsilon} = K$ and $\overline{K} = \overline{K}^\circ = K$. In particular, the scalar multiplication $\circ$ is defined over $K$ and $t\circ \lambda = t\lambda^2$ for any $t, \lambda \in K$. The additive group of $K$ equipped with $\circ$ as the scalar multiplication is a $K$-vector space. In order to distinguish between this vector space and the field $K$ itself we denote the latter by the letter $K$, keeping the symbol $\overline{K}$ for the vector space structure $(K,\circ)$. Given an element $t\in K$, if we regard it as a vector of $\overline{K}$ then we write $\bar{t}$ rather than $t$. 

Put $\widetilde{V} := V\oplus \overline{K}$. The set $W := \{x+\overline{\gamma_E(x)}\}_{[x]\in e(P)}$ is a subset of $\widetilde{V}$ and contains $E$. However $E$ spans $V$, the latter being now regarded as a subspace of $\widetilde{V}$. Therefore $\langle W\rangle\supseteq V$. It follows that $\langle W\rangle$ also contains the set $\{\overline{\gamma_E(x)}\}_{[x]\in e(P)}$. The latter spans $\overline{R}$ and $\overline{R} = \overline{K}$, by assumption. Therefore $W$ spans $\widetilde{V}$. We now define a quadratic form $\tilde{q}$ and an alternating form $\tilde{f}$ on $\widetilde{V}$, as follows:
\[\tilde{q}(x+\bar{t}) = \gamma_E(x) + t ~~ \mbox{for any}~ x\in V ~\mbox{and}~ \bar{t}\in\overline{K}.\]
\[\tilde{f}(x+\bar{t}, y+\bar{s}) = f(x,y) ~~ \mbox{for any}~ x, y \in V ~\mbox{and}~ \bar{t}, \bar{s}\in K.\]
It is readily seen that $\tilde{q}$ is indeed a quadratic form and $\tilde{f}$ is its sesquilinearization. Note that $\mathrm{Rad}(f) = \overline{K}$ and $\overline{K}$ contains no $\tilde{q}$-singular point. Hence $\tilde{q}$ is non-singular. Accordingly, the polar space $S_{\tilde{q}} = (P_{\tilde{q}},L_{\tilde{q}})$ associated to $\tilde{q}$ in $\mathrm{PG}(\widetilde{V})$ is non-degenerate. Moreover $S_{\tilde{q}}$ is a subspace of the polar space  $S_{\tilde{f}}$ associated to $\tilde{f}$, as $\tilde{f}$ is the sesquilinearization of $\tilde{q}$.  

For $x\in V$ and $\bar{t}\in \overline{K}$ we have $\tilde{q}(x+\bar{t}) = 0$ if and only if $t = \gamma_E(x)$. Hence the set $\widetilde{P} := \{[v]\}_{v\in W}$ is contained in $P_{\tilde{q}}$. It is not difficult to see that $\widetilde{P}$ is a subspace of $S_{\tilde{q}}$. Let $\widetilde{S}$ be the polar space induced by $S_{\tilde{q}}$ on $\widetilde{P}$. Clearly, $\widetilde{S}$ is a subspace of $S_{\tilde{q}}$. Hence it is also a subspace of $S_{\tilde{f}}$, since $S_{\tilde{q}}$ is a subspace of $S_{\tilde{f}}$. Since $\widetilde{P}$ spans $\widetilde{V}$ and $\widetilde{S}$ is a subspace of $S_{\tilde{q}}$, the radical of $\widetilde{S}$ is contained in the radical of $S_{\tilde{q}}$. However $S_{\tilde{q}}$ is non-degenerate. Hence $\widetilde{S}$ is non-degenerate. Consequently, property (E1) (whence (E2), (E3) and (E4)) hold for $\widetilde{S}$ and $\tilde{f}$. 

We shall prove that $\widetilde{S} = S_{\tilde{q}}$. By way of contratiction, let $[a]\in P_{\tilde{q}}-\widetilde{P}$. Note that $a\not\in \mathrm{Rad}(\tilde{f})$, because $S_{\tilde{q}}$ is non-degenerate. Then, by $(E4)$ applied to $\widetilde{S}$ and $\tilde{f}$, there is a line $l$ of $\mathrm{PG}(\widetilde{V})$ containing $[a]$ and two distinct points $[b], [c]\in\widetilde{P}$. The line $l$ belongs to $L_{\tilde{q}}$, since it contains at least three distinct points of $P_{\tilde{q}}$ and $\tilde{q}$ is quadratic. Consequently, $l$ is totally singular for $\tilde{q}$. Hence $l$ is also totally isotropic for $\tilde{f}$. In particular $f(b,c) = 0$. This forces $l$ to be a line of $\widetilde{S}$ too, a contradiction with the choice of $[a]\not\in \widetilde{P}$. Therefore $\widetilde{S} = S_{\tilde{q}}$.

The projection $\pi_{\overline{K}}:\widetilde{V}\rightarrow \widetilde{V}/\overline{K} = V$ induces an isomorphism from $\widetilde{S}$ to $e(S)$. On the other hand, the quotient $\tilde{q}_{\overline{K}}$ of $\tilde{q}$ by $\overline{K}$ is trivial. Hence $\pi_{\overline{K}}$ induces an isomorphism from $S_{\tilde{q}}$ to $S_f$, by claim (2) of Theorem \ref{quot6}. However $S_{\tilde{q}} = \widetilde{S}$. Therefore $e(S) = S_f$.  \eop 

\section{Initial embeddings} 

In this section we shall revisit Theorem \ref{Theo 3}, giving an elementary proof the fact that the embeddings considered in Theorem \ref{Theo 3} are dominant and a proof of the last claim of Theorem \ref{Theo 3} in the case of rank at least $3$, different from the original proof of Tits \cite{Tits}. 

With $e:S\rightarrow \mathrm{PG}(V)$ and $f:V\times V\rightarrow K$ as in the previous section, let $q:V\rightarrow \overline{K}/\overline{R}$ be the generalized pesudo-quadratic form defined as in (\ref{definition q}). By Theorems \ref{main1} and \ref{qE4}, either $q$ is non-trivial and $e(S) = S_q$ or $K$ is a field, $f$ is alternating and $e(S) = S_f$. 

The existence of the cover $q^{\overline{R},\{\bar{0}\}}$ makes it clear that, if $e(S) = S_q$, then $e$ is dominant only if $\overline{R} = \{\bar{0}\}$, namely $q$ is pseudo-quadratic. Conversely, 

\begin{lemma}\label{main2}
Suppose that either $q$ is pseudo-quadratic or $f$ is alternating and $\mathrm{char}(K)\neq 2$. Then $e$ is dominant.
\end{lemma}
\pr This lemma is contained in Theorem \ref{Theo 3} but, since we are revisiting Theorem \ref{Theo 3}, we shall give a proof independent of that theorem.  Our proof exploits Theorems \ref{main1} and \ref{qE4} and properties of quotients of generalized pseudo-quadratic forms. 

Let $\tilde{e}:S\rightarrow\mathrm{PG}(\widetilde{V})$ be the hull of $e$. Then there exists a reflexive sesquilinear form $\tilde{f}:\widetilde{V}\times\widetilde{V}\rightarrow K$ such that $\tilde{e}(S)$ is a subspace of $S_{\tilde{f}}$. Let $\tilde{q}:\widetilde{V}\rightarrow \overline{K}/\overline{R}$ be the generalized pseudo-quadratic form defined as in (\ref{definition q}) but with $V$ and $f$ replaced with $\widetilde{V}$ and $\tilde{f}$ respectively. By Theorems \ref{main1} and \ref{qE4}, either $\overline{R} \subset \overline{K}$ and $\tilde{e}(S) = S_{\tilde{q}}$ or $\overline{R} = \overline{K}$ and $\tilde{e}(S) = S_{\tilde{f}}$.  

As $\tilde{e}$ is the hull of $e$, there exists a subspace $U$ of $\mathrm{Rad}(\tilde{f})$ such that $e \cong \tilde{e}/U$. If $\tilde{e}(S) = S_{\tilde{f}}$ then $\tilde{f}$ is non-degenerate. In this case $U = \{0\}$, whence $e \cong \tilde{e}$, namely $e$ is dominant. 

Suppose that $\overline{R} \subset \overline{K}$. Then $\tilde{e}(S) = S_{\tilde{q}}$ and $e(S) = S_{\tilde{q}_U}$, where $\tilde{q}_U$ is the quotient of $\tilde{q}$ by $U$, regarded as a generalized pseudo-quadratic form on $V$ via an obvious identification of $V$ with $\widetilde{V}/U$. Then $q$ and $\tilde{q}_U$ are proportional, by Theorem \ref{proportional gen quadric theo}. Hence $\tilde{q}_U$ is a pseudo-quadratic form. However pseudo-quadratic forms do not admit proper covers (Subsection \ref{dominant covers}), while $\tilde{q}$ is a cover of $\tilde{q}_U$ by Theorem \ref{covers-main-ter}. Hence $\tilde{q}_U \cong \tilde{q}$, namely $U = \{0\}$. Again, $e \cong \tilde{e}$.  \eop

\bigskip 

Turning back to the general case, when $\overline{R}\subset \overline{K}$ we denote by $\tilde{e}$ the composition of $e$ with the lifting of $e_q:S_q\rightarrow \mathrm{PG}(V)$ to $\widetilde{V} := V^{\overline{R}}$. Thus, $\tilde{e}(S) = S_{\tilde{q}}$, where $\tilde{q} := q^{\overline{R},\{\bar{0}\}}$ is the dominant cover of $q$, as in the proof of Theorem \ref{main1}. When $\overline{R} = \overline{K}$ and $\mathrm{char}(K)\neq 2$ we set $\widetilde{V} := V$ and $\tilde{e} := e$. Finally, let $\overline{R} = \overline{K}$ but $\mathrm{char}(K) = 2$. It is well known that in this case $e$ is a quotient of an embedding $\tilde{e}:S\rightarrow \mathrm{PG}(\tilde{V})$, where $\tilde{e}(S) = S_{\tilde{q}}$ for a suitable quadratic form $\tilde{q}:\widetilde{V}\rightarrow K$ (see e.g. De Bruyn and Pasini \cite{DBP}).

\begin{theo}\label{main3} 
In each of the cases considered above the embedding $\tilde{e}$ is dominant, whence it is the hull of $e$.  
\end{theo}
\pr This statement immediately follows from Lemma \ref{main2}, recalling that dominat covers of generalized pseudo-quadratic forms are pseudo-quadratic forms (Subsection \ref{dominant covers}). \eop   

\begin{co}\label{main4}
With $\tilde{e}$ as above, assume moreover that $S$ has rank at least $3$. Then $\tilde{e}$ is absolutely initial.
\end{co}
\pr Embeddable polar spaces of rank at least 3 satisfy the conditions of the main theorem of Kasikova and Shult \cite{KS}, which are sufficient for the existence of a $K$-initial embedding. Therefore the embedding $\tilde{e}$, being dominant, is also also $K$-initial (see Subsection \ref{Emb}). On the other hand, since $K$ is the divison ring coordinatizing the planes of $S$, all projective embeddings of $S$ are $K$-embeddings, namely $S$ is defined over $K$. Hence $\tilde{e}$ is absolutely initial.  \eop 

 \bigskip

The statement of Corollary \ref{main4} is included in Theorem \ref{Theo 3}, which is a rephrasing of Theorem 8.6 of Tits \cite{Tits}, but the proof given by Tits for Theorem 8.6 of \cite{Tits} is rather different from our proof of Corollary \ref{main4}. In our proof we rely on the main result of Kasikova and Shult \cite{KS}, which can be applied to polar spaces of rank at least 3 thanks to the fact that the maximal singular subspaces of such a polar space are projective spaces of dimension at least 2, while the proof given by Tits in \cite{Tits} relies on certain deep properties of projective lines. Tits's proof also applies to polar spaces of rank 2 but for the two exceptional cases described in the following theorem (and mentioned in Theorem \ref{Theo 3}).  

\begin{theo}\label{main5}{\rm [Tits \cite[8.6]{Tits}]}
The embedding $\tilde{e}$ is absolutely initial even if $S$ has rank $2$, except in the following two cases: 

\begin{itemize}
\item[$(1)$] $S$ is a grid and $|K| > 4$.
\item[$(2)$] $K$ is a quaternion division ring, $\widetilde{V} = V(4,K)$ and, modulo proportionality and isomorphisms, $\varepsilon = -1$, $\sigma$ is the standard involution of $K$, we have $K_{\sigma,\varepsilon} = Z(K)$ and $\tilde{q}(x_1,x_2,x_3,x_4) = x_1^\sigma x_2+ x_3^\sigma x_4 + K_{\sigma,\varepsilon}$ for every vector $(x_1,x_2,x_3,x_4)\in \widetilde{V}$. 
\end{itemize}
In case $(1)$ we have as many isomorphism classes of projective embeddings as the cosets of $\mathrm{P}\Gamma\mathrm{L}(2,K)$ in the group of all permutations of the set $\mathrm{PG}(1,K)$. In case $(2)$ only two isomorphism classes of projective embeddings exist. In either case, all projective embeddings of $S$ are dominant. 
\end{theo}


\bigskip

\noindent
address of the author\\

\noindent
Antonio Pasini,\\
Department of Information Engeering and Mathematics,\\
Uiversity of Siena\\
Via Roma 56, 53100 Siena, Italy\\
antonio.pasini@unisi.it 

\begin{thebibliography}{XXX}
\bibitem{BC} F. Buekenhout and A. M. Cohen. {\em Diagram Geometries}, Springer, Berlin, 2013. 
\bibitem{BL} F. Buekenhout and C. Lef\'{e}vre. Generalized quadrangles in projective spaces. {\em Arch. Math.} {\bf 25} (1974), 540-552. 
\bibitem{BS} F. Buekenhout and E. E. Shult. On the foundations of polar geometry. {\em Geometriae Dedicata} {\bf 3} (1974), 155-170. 
\bibitem{DBP} B. De Bruyn and A. Pasini. On symplectic polar spaces over non-perfect fields of characteristic $2$. {\em Linear and Multilinear Algebra} {\bf 47} (2009), 567-575. 
\bibitem{KS} A. Kasikova and E. E. Shult. Absolute embeddings of point-line geometries. {\em J. Algebra} {\bf 238} (2001), 100-117. 
\bibitem{Ron} M. A. Ronan. Embeddings and hyperplanes of discrete geometries. {\em European J. Combin.} {\bf 8} (1987), 179-185.
\bibitem{Tits} J. Tits. {\em Buildings of Shperical Type and Finite $BN$-pairs}. Springer Lecture Notes {\bf 386} (1974), Springer, Berlin. 
\bibitem{Veldkamp} F. D. Veldkamp. Polar Geometry I-V. {\em Indag. Math.} {\bf 21} (1959), 512-551 and {\bf 22} (1959), 207-212.  
\end{thebibliography}
\end{document}